\documentclass[11pt, a4paper]{article}

\usepackage[utf8]{inputenc}
\usepackage[T1]{fontenc}
\usepackage{amssymb, amsthm, amsfonts, mathrsfs}
\usepackage{geometry}
\usepackage{mathtools}
\usepackage[all]{xy}
\usepackage{graphicx}
\usepackage{microtype}
\usepackage{enumitem}
\usepackage[hidelinks]{hyperref}

\theoremstyle{plain}
\newtheorem{theorem}{Theorem}[section]
\newtheorem{lemma}[theorem]{Lemma}
\newtheorem{proposition}[theorem]{Proposition}
\newtheorem{corollary}[theorem]{Corollary}
\newtheorem*{proposition*}{Proposition}
\theoremstyle{definition}
\newtheorem{definition}[theorem]{Definition}
\newtheorem{example}[theorem]{Example}

\theoremstyle{remark}
\newtheorem{remark}[theorem]{Remark}
\newcommand{\Z}{\mathbb{Z}}
\newcommand{\R}{\mathbb{R}}
\newcommand{\C}{\mathbb{C}}
\newcommand{\N}{\mathbb{N}}
\newcommand{\RP}{\mathbb{RP}}
\newcommand{\Oisz}{\mathcal{O}^{isz}(\C)}
\newcommand{\Ostar}{\mathcal{O}^*(\C)}
\newcommand{\Clf}{C^{lf}_{\infty}(\C)}
\newcommand{\Conf}{Conf^{lf}_{\infty}(\C)}

\newcommand{\cl}[1]{\overline{#1}}

\DeclareMathOperator{\Aut}{Aut}
\DeclareMathOperator{\supp}{supp}
\DeclareMathOperator{\dist}{dist}
\DeclareMathOperator{\id}{id}
\DeclareMathOperator{\Hom}{Hom}

\title{The asphericity of locally finite infinite configuration spaces and Weierstrass entire coverings}

\author{Jyh-Haur Teh\thanks{Department of Mathematics, National Tsing Hua University, Taiwan.
E-mail: \texttt{jyhhaur@math.nthu.edu.tw}.}}
\date{} % leave empty

\begin{document}
%\date{Department of Mathematics, National Tsing Hua University of Taiwan\\
%Email: \ jyhhaur@math.nthu.edu.tw\\}
\maketitle

\begin{abstract}
Let $\Conf$ and $\Clf$ denote the locally finite infinite ordered and unordered
configuration spaces of the complex plane. We prove that both $\Conf$ and $\Clf$ are aspherical.
We further obtain a locally finite analogue of the braid exact sequence,
\[
1\longrightarrow H^{lf}(\infty)\longrightarrow B^{lf}(\infty)\longrightarrow \Aut(\N)\longrightarrow 1,
\]
where $H^{lf}(\infty)=\pi_1(\Conf)$ and $B^{lf}(\infty)=\pi_1(\Conf//\Aut(\N))$, the fundamental group of the homotopy
quotient of $\Conf$ by $\Aut(\N)$.
Building on this, we classify connected countably infinite--sheeted covering spaces and give a criterion for when such
a covering can be realized from the zero set of a family of entire functions $F:X\times\C\to\C$.
In particular, if $\pi_1(X)$ is free and $H^2(X;\Z)=0$, then every countably infinite--sheeted covering space
over $X$ is realizable.
\end{abstract}

\section{Introduction}
Configuration spaces of points on a space $M$ form a classical meeting ground for geometry, topology, and
analysis. For $n\ge 1$, the ordered configuration space
\[
Conf_n(M)=\{(x_1,\dots,x_n)\in M^n \mid x_i\neq x_j\text{ for }i\neq j\}
\]
and its unordered quotient $C_n(M)=Conf_n(M)/S_n$ encode, in particular, the (pure) braid
groups of $M$ when $M$ is a surface. The foundational fibrations of Fadell--Neuwirth \cite{FadellNeuwirth1962},
together with later developments (e.g.\ \cite{McDuff1975, Segal1973}), make these spaces central objects in
algebraic topology. For background and many further directions we refer to Kallel's survey
\cite{KallelUserGuide} and the references therein; see also \cite{KallelGTM08} for a geometric approach to
braid spaces and their homology.

In this paper we study an ``infinite, locally finite'' analogue of planar configuration spaces.
Write $\Conf$ for the space of \emph{locally finite ordered} configurations
$\widetilde{A}=(a_1,a_2,\dots)$ in $\C$ (no repetitions and no accumulation in compact subsets). Let $\Clf$
be the corresponding space
of \emph{unordered} locally finite configurations (i.e.\ infinite locally finite subsets of $\C$), obtained as the
quotient $\Conf/\Aut(\N)$. A guiding principle is that, unlike the standard direct limit
$\varinjlim_nConf_n(\C)$, the spaces $\Conf$ and $\Clf$ capture a genuine ``configuration of
countably many particles escaping to infinity'' and carry a topology compatible with the \emph{vague
topology} obtained by viewing configurations as counting measures. Moreover, $\pi_1\Conf$ is
uncountable, whereas $\pi_1(\varinjlim_nConf_n(\C))$ is countable.

This project originated in our search for a version of Lawson homology and morphic cohomology
\cite{Lawson89,LawsonSDG93,Friedlander91,FriedlanderLawson92,FriedlanderLawson97,Teh08TransAMS,TehYang20CM,TehYang21CM}
for analytic cycles on unbounded analytic varieties. This led us to a basic but, to our knowledge, is not explicitly
addressed in the literature: what is the topology of the space of complex entire functions with infinitely many simple zeros?
We were not able to locate a reference treating this problem \cite{Birman74, KasselTuraev08, FadellHusseini01}. Informal discussions with specialists in braid theory suggest that the question is not currently treated in the literature in the form needed for our applications.
The results of this paper provide a solution in the planar case; the necessary constructions are somewhat delicate, but they
also indicate a route to higher-dimensional analogues.

Beyond its intrinsic geometric interest, the \emph{unordered locally finite configuration space}
\[
\Gamma(\C)\;=\;\Bigl\{\gamma\subset \C \ \Big|\ \#(\gamma\cap K)<\infty \text{ for every compact }K\Subset\C\Bigr\}
\]
appears in several other areas of mathematics. Identifying $\gamma$ with the \emph{simple counting measure}
$\mu_\gamma=\sum_{x\in\gamma}\delta_x$, one may topologize $\Gamma(\C)$ by the \emph{vague topology} (weak
convergence against test functions in $C_c(\C)$). This is the standard state space in point process theory
and stochastic geometry \cite{KallenbergRM,DaleyVereJonesII}; in particular, $\Gamma(\C)$ serves as the
canonical sample space for Poisson and more structured processes \cite{DaleyVereJonesII,KallenbergRM}. In
random matrix theory, the eigenvalues of complex Ginibre matrices form a determinantal point process in the
plane \cite{Ginibre1965,ForresterLogGases}. In random analytic function theory, zeros of Gaussian analytic
functions and related ensembles give natural random configurations; they are closely connected to
determinantal point processes and provide a probabilistic complex-analytic way to build measures on
$\Gamma(\C)$ \cite{HKPVBook,KabluchkoZaporozhets2014}.

In continuum statistical mechanics, a ``gas'' of indistinguishable particles in the plane is described by
probability measures on $\Gamma(\C)$. The foundational framework for Gibbs measures and phase transitions is
typically formulated on configuration spaces \cite{Ruelle1969,Georgii2011}. Beyond equilibrium, one studies
Markovian dynamics on $\Gamma(\C)$ (Kawasaki dynamics, Glauber-type dynamics, interacting diffusions, \dots).
A large body of work develops analysis and geometry on configuration spaces (Dirichlet forms,
integration-by-parts formulas, differential structures) precisely in order to construct and analyze such
dynamics \cite{AlbeverioKondratievRoeckner1998}. In the planar/logarithmic interaction setting relevant to
random matrices, infinite-dimensional SDEs whose unlabeled state lives in $\Gamma(\C)$ have been constructed
and studied by Osada and others \cite{Osada2012}.

For finite $n$, the ordered configuration space $Conf_n(\C)$ is a classical model for braid groups and has
rich cohomology \cite{FadellNeuwirth1962,Arnold1969}. Configuration spaces also arise naturally in
collision-free motion planning: the configuration space of distinct points encodes feasible multi-agent
trajectories, with braids recording entanglement/avoidance patterns \cite{Ghrist2010}. For a broad
perspective on configuration spaces, including many further applications, see Kallel's user guide
\cite{KallelUserGuide2024}.

\medskip
\noindent\textbf{Asphericity.}
Our first main result is that both $\Clf$ and $\Conf$ are aspherical; in particular, they are
Eilenberg--Mac~Lane spaces of type $K(\pi,1)$. For $\C$ this parallels the classical fact that
$Conf_n(\C)$ is aspherical for each finite $n$ (pure braid groups), but the locally finite infinite
setting requires controlling infinitely many points simultaneously. The key input is an
\emph{Infinite Gap Lemma}: any map $S^k\to C^{lf}_\infty(\C)$ ($k\ge2$) is homotopic,
via compactly supported ambient isotopies in pairwise disjoint annuli, to a map whose configurations
avoid infinitely many prescribed radial annuli.
This produces a decomposition into disjoint shells on which the map factors through finite configuration spaces
in a disk/annulus; using that these finite configuration spaces are $K(\pi,1)$, we null-homotope shell-by-shell
and concatenate the homotopies, with vague continuity guaranteed because each test function has compact support
and meets only finitely many shells. The same shellwise lifting and concatenation, using $\pi_1(S^k)=0$ to lift
to ordered configurations, yields the asphericity of $Conf^{lf}_\infty(\C)$ as well.

\medskip
\noindent\textbf{Zero sets of entire functions.}
A second theme is the relationship with classical complex analysis. Let $\mathcal O^{\mathrm{isz}}(\C)$ be the
space of entire functions with infinitely many simple zeros, and let $\mathcal O^{*}(\C)$ be
the multiplicative group of nowhere--vanishing entire functions. In \S2 we show that the zero--set map
induces a homeomorphism between the quotient
$\mathcal O^{\mathrm{isz}}(\C)/\mathcal O^{*}(\C)$ and $\Clf$. The inverse is constructed via Weierstrass
products.

\medskip
\noindent\textbf{Countable coverings and Weierstrass realizability.}
The last part of the paper concerns the classification and realization of countably infinite--sheeted
coverings. Given a connected, locally path--connected, semilocally simply connected space $X$, isomorphism
classes of connected countable coverings of $X$ correspond to conjugacy classes of homomorphisms
$\pi_1(X)\to\Aut(\N)$. We introduce a class of \emph{embedded} countable coverings and relate them to maps
into the Borel construction $\Conf//\Aut(\N)$, yielding a concrete ``universal'' parameter space for such
coverings. We then ask when an embedded covering can be realized as the zero locus $Z(F)\subset X\times\C$
of a \emph{Weierstrass entire family} $F:X\times\C\to\C$. Our main criterion
(Theorem~\ref{thm:criterion}) identifies an obstruction class $c_1(E)\in H^2(X;\Z)$ whose vanishing is
equivalent to the existence of a global Weierstrass realization. In particular, when $H^2(X;\Z)=0$, every
closed embedded countable covering is Weierstrass--realizable.

\medskip
\noindent\textbf{Organization of the paper.}
In Sections~2--3 we construct the zero--set homeomorphism and establish its continuity properties,
thereby confirming that the vague topology on \(C^{lf}_{\infty}(\C)\) is the appropriate one.
We then introduce a natural metric on \(\Conf\), the sum metric, under which the projection
\(P:\Conf\to\Clf\) becomes continuous. We also show that \(\pi_1(\Conf)\) is uncountable and
for each $n\in \N$, \(\pi_1(Conf_n(\C))\) embeds into \(\pi_1(\Conf)\).

Section~4 proves path-connectedness of \(\Conf\) and \(\Clf\) via iterated isotopies.
In Section~5 we establish the infinite gap lemma and use it to prove that \(\Clf\) is aspherical.
Section~6 develops the notions of ordered gaps and shell decompositions to prove asphericity of \(\Conf\).

In Section~7 we study countably infinite-sheeted coverings and replace \(P:\Conf\to\Clf\) by a
homotopy-quotient map \(\overline{P}:\Conf//\Aut(\N)\to\Clf\).
Section~8 introduces Weierstrass entire coverings and resolves the realization problem.
An appendix shows that, for finite unordered configuration spaces, the usual Euclidean topology
agrees with the vague topology.

\section*{Acknowledgements}
The author thanks Professor \textsc{Sadok Kallel} for his interest in this work and for helpful
correspondence.
His writings on configuration spaces, in particular the survey \cite{KallelUserGuide} and related work
\cite{KallelGTM08, KallelTaamallah13}, provided valuable context and motivation.

\section{Zero set homeomorphism}

Let $\Oisz$ denote the collection of entire functions on $\C$ possessing infinitely many simple zeros. Let $\Ostar$ denote the collection of non-vanishing entire functions. The space of entire functions is endowed with the standard topology of uniform convergence on compact subsets.

We consider the quotient space $\Oisz / \Ostar$
equipped with the quotient topology induced by the uniform metric $d_U$ on compacta.

Let $\Clf$ be the collection of countably infinite, locally finite subsets of $\C$. We identify a set $A \in \Clf$ with the positive divisor (Radon measure) $\mu_A = \sum_{a \in A} \delta_a$. The appropriate topology for this space is the \textbf{vague topology} (weak-* convergence of measures).

For $m\in\N$ let $K_m:=\{z\in\C:\ |z|\le m\}$. Denote by $C_c(\C)=\{\psi : \psi: \C \rightarrow \R \mbox{ is continuous with compact support}\}$.

\begin{lemma}\label{lem:testfamily}
There exists a sequence $(\varphi_j)_{j\ge1}\subset C_c(\C)$ such that for every $m\in\N$,
the subfamily $\{\varphi_j:\supp(\varphi_j)\subset K_m\}$ is dense in
\[
C_c(K_m):=\{\psi\in C_c(\C):\supp(\psi)\subset K_m\}
\]
with respect to the global supremum norm $\|\cdot\|_\infty$.
\end{lemma}

\begin{proof}
Fix $m$ and write $\partial K_m=\{z\in\C:|z|=m\}$.
If $\psi\in C_c(K_m)$, then $\psi\equiv0$ on $\C\setminus K_m$; hence by continuity we must have
$\psi|_{\partial K_m}\equiv 0$. Therefore the restriction maps $C_c(K_m)$
isometrically into the closed subspace
\[
C_{\partial}(K_m):=\{f\in C(K_m): f|_{\partial K_m}\equiv 0\}
\subset C(K_m).
\]
Conversely, if $f\in C_{\partial}(K_m)$, then the extension $\tilde f$ defined by
$\tilde f|_{K_m}=f$ and $\tilde f\equiv 0$ on $\C\setminus K_m$ is continuous on $\C$,
hence $\tilde f\in C_c(\C)$ with $\supp(\tilde f)\subset K_m$. We have $C_c(K_m)=C_{\partial}(K_m)$.

Since $C_{\partial}(K_m)$ is a closed subspace of the separable Banach space $C(K_m)$, it is separable.
Choose a countable dense subset $\{f_{m,k}\}_{k\ge1}\subset C_{\partial}(K_m)$.
Extend each $f_{m,k}$ by $0$ outside $K_m$ to obtain $\varphi_{m,k}\in C_c(\C)$ with $\supp(\varphi_{m,k})\subset K_m$.

Finally, enumerate the countable set $\{\varphi_{m,k}:m,k\in\N\}$ as a single sequence $(\varphi_j)_{j\ge1}$.
Density in each $C_c(K_m)$ follows from the preceding identification.
\end{proof}

\begin{definition}[The Vague Metric $d_{\mathcal{V}}$]
Throughout this paper, we fix $\{\varphi_j\}_{j=1}^{\infty}$ to be the countable subset of $C_c(\C)$ that is dense in $C_c(\C)$ with respect to the supremum norm as constructed above.

We define the vague metric $d_{\mathcal{V}}$ on $\Clf$ by:
\[
d_{\mathcal{V}}(A, B) = \sum_{j=1}^{\infty} \frac{1}{2^j} \frac{\left| \sum_{a \in A} \varphi_j(a) - \sum_{b \in B} \varphi_j(b) \right|}{1 + \left| \sum_{a \in A} \varphi_j(a) - \sum_{b \in B} \varphi_j(b) \right|}
\]
\end{definition}

\begin{proposition}\label{prop:dVmetric}
$d_{\mathcal V}$ is a metric on $\Clf$.

\end{proposition}

\begin{proof}
Probably the only nontrivial part is to show that if $A, B\in \Clf$ satisfy $d_{\mathcal V}(A,B)=0$, then $A=B$.
Assume $d_{\mathcal V}(A,B)=0$. Then for each $j$,
\[
\sum_{a\in A}\varphi_j(a)=\sum_{b\in B}\varphi_j(b).
\]
Let $\psi\in C_c(\C)$ be arbitrary and choose $m$ with $\supp(\psi)\subset K_m$.
By Lemma~\ref{lem:testfamily} there exists a sequence $(\varphi_{j_n})_{n\geq 1}$ with $\supp(\varphi_{j_n})\subset K_m$ and
$\|\varphi_{j_n}-\psi\|_\infty\to0$.
Since $A\cap K_m$ and $B\cap K_m$ are finite, we may pass to the limit:
\[
\sum_{a\in A}\psi(a)=\lim_{n\to\infty}\sum_{a\in A}\varphi_{j_n}(a)
=\lim_{n\to\infty}\sum_{b\in B}\varphi_{j_n}(b)=\sum_{b\in B}\psi(b).
\]
Thus $\int \psi\,d\mu_A=\int \psi\,d\mu_B$ for all $\psi\in C_c(\C)$, hence $\mu_A=\mu_B$ as Radon measures.

To show $A=B$, fix $a\in A$. Since $A$ is locally finite, there exists $r_0>0$ such that
$A\cap \cl{B(a,r_0)}=\{a\}$, hence $\mu_A(B(a,r))=1$ for all $0<r\le r_0$.
Since $\mu_A=\mu_B$, we also have $\mu_B(B(a,r))=1$ for all $0<r\le r_0$.
In particular, for every $n$ with $1/n\le r_0$ the ball $B(a,1/n)$ contains a point of $B$.
Therefore $a\in \cl{B}$. Since the set $B$ is closed, so $a\in B$.
Hence $A\subset B$, and by symmetry $B\subset A$. Thus $A=B$.
\end{proof}

\begin{definition}[The Root Map]
The root map $Z: \Oisz / \Ostar \to (\Clf, d_{\mathcal{V}})$ is defined by
\[
Z([f]) = f^{-1}(0)
\]
\end{definition}

\begin{theorem}
The root map $Z: \Oisz / \Ostar \to (\Clf, d_{\mathcal{V}})$ is a homeomorphism.
\end{theorem}

\begin{proof}
It is clear that $Z$ is well defined.
We proceed by showing that $Z$ is a bijection, and that both $Z$ and $Z^{-1}$ are sequentially continuous.

\paragraph{Surjectivity:}
Let $A \in \Clf$. Since $A$ is a discrete set in $\C$, the Weierstrass Product Theorem guarantees the existence of an entire function $f$ vanishing exactly on $A$ with multiplicity $1$. Since $A$ is infinite, $f \in \Oisz$. Thus $Z([f]) = A$.

\paragraph{Injectivity:}
Suppose $Z([f]) = Z([g])$. Then $f$ and $g$ have identical simple zeros. Define $h(z) = f(z)/g(z)$. The singularities of $h$ are removable, and the resulting extension is non-vanishing on $\C$. Thus $h \in \Ostar$, implying $f = g \cdot h$, so $[f] = [g]$.

\paragraph{Continuity of $Z$:}

Let $\{[f_k]\}_{k \in \N}$ be a sequence in $\Oisz / \Ostar$ such that $[f_k] \to [f]$.
Since the $\mathcal{O}^*(\C)$-action on $\mathcal{O}^{isz}(\C)$ is by homeomorphisms, the quotient map
$\mathcal{O}^{isz}(\C) \longrightarrow \mathcal{O}^{isz}(\C)/\mathcal{O}^*(\C)$
is open. After choosing a representative $f\in [f]$, one can choose representatives $f_k\in [f_k]$ such that $f_k\to f$ uniformly on compact subsets of $\mathbb{C}$.

Let $\varphi \in C_c(\C)$ be a test function from the definition of the metric $d_{\mathcal{V}}$. Let $K = \operatorname{supp}(\varphi)$. Choose a bounded open set $G$ such that $K \subset G$ and $\partial G \cap Z(f) = \emptyset$.

Since $f_k \to f$ uniformly on the compact set $\partial G$ and $\min_{z \in \partial G} |f(z)| > 0$, Hurwitz's Theorem implies that for sufficiently large $k$:
\begin{itemize}
    \item The number of zeros of $f_k$ in $G$ equals the number of zeros of $f$ in $G$.
    \item The zeros of $f_k$ in $G$ converge pointwise to the zeros of $f$ in $G$.
\end{itemize}
Consequently, since $\varphi$ is continuous:
\[
\lim_{k \to \infty} \sum_{a \in Z(f_k)} \varphi(a) = \sum_{a \in Z(f)} \varphi(a).
\]
Since this holds for every $\varphi_j$ in the definition of the metric, $d_{\mathcal{V}}(Z(f_k), Z(f)) \to 0$.

\paragraph{Continuity of $Z^{-1}$:}
Let $A_k, A\in \Clf$. Assume $d_{\mathcal{V}}(A_k, A) \to 0$. We must construct representatives $f_k, f$ such that $Z(f_k) = A_k$, $Z(f) = A$, and $f_k \to f$ uniformly on compact subsets of $\C$.

\begin{description}
\item[Step 1: Indexing:]
Fix an order of elements of $A$. Write $A = \{a_n\}_{n=1}^{\infty}$. For any disk $D_R$ with boundary disjoint from $A$, the points of $A_k \cap D_R$ eventually correspond one-to-one with $A \cap D_R$. Thus, for sufficiently large $k$, we can enumerate $A_k = \{a_{n,k}\}_{n=1}^{\infty}$ such that:
\[
\lim_{k \to \infty} a_{n,k} = a_n \quad \text{for every } n \in \N.
\]

\item[Step 2: Constructing $f$:]
Define the factors $h_n(z)$:
\[
h_n(z) =
\begin{cases}
z & \text{if } a_n = 0, \\
E_n\left(\frac{z}{a_n}\right) & \text{if } a_n \neq 0.
\end{cases}
\]
where $E_n\left(\frac{z}{a_n}\right) = \left(1 - \frac{z}{a_n}\right) \exp\left( \sum_{j=1}^n \frac{1}{j} (\frac{z}{a_n})^j \right)$
is the Weierstrass elementary factor.
Set $f(z) = \prod_{n=1}^{\infty} h_n(z)$.

\item[Step 3: Constructing $f_k$:]
For $k$ large enough such that $a_{n,k} \neq 0$ whenever $a_n \neq 0$, define $f_k$ term-by-term to match the structure of $f$:
\[
h_{n,k}(z) =
\begin{cases}
(z - a_{n,k}) & \text{if } a_n = 0, \\
E_n\left(\frac{z}{a_{n,k}}\right) & \text{if } a_n \neq 0.
\end{cases}
\]
Set $f_k(z) = \prod_{n=1}^{\infty} h_{n,k}(z)$.

\item[Step 4: Uniform Convergence:]
Let $K \subset \C$ be compact. Choose $R > 0$ such that $K \subset \{|z| < R\}$. Since $|a_n| \to \infty$, choose $N$ such that $|a_n| > 2R$ for all $n > N$, and assume without loss of generality that $\partial D_{2R} \cap A = \emptyset$.

\paragraph{Conservation of Number (Tail Control):}
Since $A_k \to A$ vaguely, $\lim_{k \to \infty} \#(A_k \cap D_{2R}) = \#(A \cap D_{2R}) = N$.
Since we have indexed the first $N$ points such that $a_{n,k} \to a_n$ (which lie inside $D_{2R}$), these account for all zeros in the disk. Therefore, for all $n > N$, the remaining zeros must lie outside:
\[
|a_{n,k}| > 2R \quad \text{for all } n > N \text{ and sufficiently large } k.
\]

\paragraph{Convergence Analysis:}
We split the product:
\[
|f_k(z) - f(z)| \le \left| \prod_{n=1}^N h_{n,k}(z) \prod_{n=N+1}^\infty h_{n,k}(z) - \prod_{n=1}^N h_n(z) \prod_{n=N+1}^\infty h_n(z) \right|
\]

\begin{enumerate}
    \item \textbf{Head ($n \le N$):} The factors depend continuously on the parameter $a$. Since $a_{n,k} \to a_n$, the finite product converges uniformly on $K$.

    \item \textbf{Tail ($n > N$):} For $z \in K$ ($|z| < R$) and $|w| > 2R$, we have $|z/w| < 1/2$. We use the standard estimate $|\log E_n(u)| \le 2|u|^{n+1}$.
    For $f_k$, using $|a_{n,k}| > 2R$:
    \[
    \sum_{n=N+1}^{\infty} |\log h_{n,k}(z)| \le \sum_{n=N+1}^{\infty} 2 \left( \frac{R}{|a_{n,k}|} \right)^{n+1} \le \sum_{n=N+1}^{\infty} 2 \left( \frac{1}{2} \right)^{n+1}<\infty
    \]
    and is independent of $k$. Since each term $h_{n,k} \to h_n$, the dominated convergence theorem for infinite products implies that the tails converge uniformly to the tail of $f$.
\end{enumerate}

Since both the head and tail converge uniformly, $f_k \to f$ uniformly on $K$. Thus $Z^{-1}$ is continuous.
\end{description}
\end{proof}

The following result shows that on finite unordered configuration spaces, the vague topology is same as the
standard Euclidean topology. From this point of view, it is coherent to use vague topology on locally finite infinite
unordered configuration space.

\begin{definition}[Finite atomic measures]
Let $U\subset\C$ be open and $n\ge 1$. Define
\[
M_n(U) := \left\{ \mu = \sum_{i=1}^n \delta_{z_i} : z_i\in U,\ z_i\neq z_j\text{ for }i\neq j \right\},
\]
viewed as a subset of the space of Radon measures on $U$, equipped with the vague topology:
$\mu_\alpha\to\mu$ if and only if $\int f\,d\mu_\alpha\to\int f\,d\mu$ for all $f\in C_c(U)$.
\end{definition}

\begin{proposition}\label{prop:finite-vague-homeo}
For any open $U\subset\C$ and $n\ge1$, there is a natural homeomorphism
\[
\overline F : C_n(U) \xrightarrow{\ \cong\ } M_n(U),
\]
where $C_n(U)$ has the quotient topology from $Conf_n(U)$ and $M_n(U)$ has the vague topology.
\end{proposition}

Since this result is not used in our main proof, we leave it to the appendix.

We will use the following standard consequence of the Fadell–Neuwirth theorem several times throughout this paper.

\begin{proposition}\label{prop:finite-aspherical}
If $U\subset\C$ is an open disc or an open annulus, then $Conf_n(U)$ is aspherical.
In particular, for $m\ge2$ and any basepoint $u_0\in Conf_n(U)$, one has $\pi_m(Conf_n(U),u_0)=0$.
\end{proposition}

\section{Continuity of the map $P$ and automorphisms}
\subsection{Continuity of the map $P: \Conf \rightarrow \Clf$}
\begin{definition}
Let $\Conf$ denote the set of all sequences $x = (x_j)_{j\in\N} \in \C^\N$
such that:
\begin{enumerate}
  \item $x_i \neq x_j$ for all $i \neq j$;
  \item (\emph{local finiteness}) for every radius $R>0$ the set
  $\{j \in \N : |x_j| \le R\}$ is finite.
\end{enumerate}
For $A=\{a_j\}^{\infty}_{j=1}\in Conf^{lf}_{\infty}(\C)$, let
$$set(A)=\{a_j\ : \ j\in \N\}$$ be the locally finite subset of $\C$ formed by terms of $A$.
There is a natural map $P: Conf^{lf}_{\infty}(\C) \rightarrow C^{lf}_{\infty}(\C)$ defined by $$P(A):=set(A)$$
\end{definition}

In other words, $\Conf$ is the space of infinite, locally finite, ordered configurations in $\C$. We always
take $\widetilde{\N}=\{j\}_{j=1}^{\infty}$ as the base point of $\Conf$. We note that this space is much more
complicated than the direct limit of finite ordered configuration spaces of $\C$.

\begin{remark}
The function
\[
d_{\mathrm{prod}}(x,y)
  \;=\; \sum_{j=1}^\infty 2^{-j}\,\min\{|x_j - y_j|,\,1\}.
\]
is a metric on $\C^{\N}$ which induces the product topology of $\C^{\N}$. But its restriction
to $\Conf$ is not the right metric as shown in the following example.
\end{remark}

\begin{example}(The sneaking zero)
Let $x = (1, 2, 3, 4, \dots)$ and $x^{(k)} = (1, 2, \dots, k, \mathbf{0}, k+2, \dots)$.
In the product topology, $x^{(k)} \to x$. However, the quotient sets are $P(x^{(k)}) = \{0, 1, \dots, k, \dots\}$ and $P(x) = \{1, 2, \dots\}$. The limit of the sets (in the vague topology) has a zero at the origin. Thus, the projection map $P$ is discontinuous.
\end{example}

\begin{definition}
Throughout this paper, we equip $\Conf$ with the sum metric
\[
d_{\sum}(x,y) = d_{\mathrm{prod}}(x,y) + d_{\mathcal V}(P(x),P(y)).
\]
\end{definition}

Convergence in $d_{\sum}$ means: coordinatewise convergence in the product metric
together with vague convergence of the underlying unlabeled configurations.
Under the sum metric, the following result is clear.

\begin{proposition}
The projection map $P: (\Conf, d_{\sum}) \to (\Clf, d_{\mathcal{V}})$ is continuous.
\end{proposition}

\subsection{Continuity of automorphisms}
\begin{definition}
Let
\[
\Aut(\N)=\{\phi:\N\to\N\mid \phi \text{ is bijective}\}.
\]
For $\phi\in\Aut(\N)$ define $\phi_*:Conf^{lf}_{\infty}(\C)\to Conf^{lf}_{\infty}(\C)$ by
\[
\phi_*((a_j)_{j\ge1})=(a_{\phi(j)})_{j\ge1}.
\]
\end{definition}

\begin{lemma}\label{lem:prodcontrol}
Fix $\phi\in\Aut(\N)$ and $\varepsilon_1>0$. Choose $N\in\N$ such that
$\sum_{j>N}2^{-j}<\frac{1}{2}\varepsilon_1$.
Let $S=\{\phi(1),\dots,\phi(N)\}$ and let $M=\max S$. Set $\delta_0=\min(\frac{1}{2}\varepsilon_1,\,1)$.
If $A=(a_j)_{j\geq 1}$ and $B=(b_j)_{j\geq 1}$ satisfy $d_{\mathrm{prod}}(A,B)<2^{-M}\delta_0$, then
\[
d_{\mathrm{prod}}(\phi_*A,\phi_*B)<\varepsilon_1.
\]
\end{lemma}

\begin{proof}
For any $k\le M$ we have
\[
2^{-k}\min(|a_k-b_k|,1)\le d_{\mathrm{prod}}(A,B)<2^{-M}\delta_0\le 2^{-k}\delta_0,
\]
hence $\min(|a_k-b_k|,1)<\delta_0$.

Now for each $1\le j\le N$ we have $\phi(j)\in S\subset\{1,\dots,M\}$, so
$\min\bigl(|a_{\phi(j)}-b_{\phi(j)}|,1\bigr)<\delta_0$.
Therefore
\begin{align*}
d_{\mathrm{prod}}(\phi_*A,\phi_*B)
&=\sum_{j=1}^{N}2^{-j}\min(|a_{\phi(j)}-b_{\phi(j)}|,1)
+\sum_{j>N}2^{-j}\min(|a_{\phi(j)}-b_{\phi(j)}|,1)\\
&<\delta_0\sum_{j=1}^{N}2^{-j}+\sum_{j>N}2^{-j}
\le \delta_0+\frac{\varepsilon_1}{2}
\le \varepsilon_1
\end{align*}

\end{proof}

\begin{theorem}\label{thm:continuity}
For every $\phi\in\Aut(\N)$, the map
\[
\phi_*:(Conf^{lf}_{\infty}(\C),d_{\sum})\longrightarrow (Conf^{lf}_{\infty}(\C),d_{\sum})
\]
is continuous. Moreover, $\phi_*$ is a homeomorphism.
\end{theorem}

\begin{proof}
Fix $\phi\in\Aut(\N)$ and let $\varepsilon>0$ be given. Put $\varepsilon_1=\varepsilon/2$.
Choose $N$ as in Lemma~\ref{lem:prodcontrol} for this $\varepsilon_1$, and let $M$ and $\delta_0$ be as there.
Define
\[
\delta:=\min\bigl(2^{-M}\delta_0,\ \varepsilon/2\bigr).
\]
Suppose $A,B\in Conf^{lf}_{\infty}(\C)$ satisfy $d_{\sum}(A,B)<\delta$. Then
\[
d_{\mathrm{prod}}(A,B)\le d_{\sum}(A,B)<2^{-M}\delta_0,
\]
so Lemma~\ref{lem:prodcontrol} implies
\begin{equation}\label{eq:prodsmall}
d_{\mathrm{prod}}(\phi_*A,\phi_*B)<\varepsilon/2.
\end{equation}
Also,
\[
d_{\mathcal V}(P(A),P(B))\le d_{\sum}(A,B)<\varepsilon/2,
\]
and
\begin{equation}\label{eq:vaguesmall}
d_{\mathcal V}(P(\phi_*A),P(\phi_*B))=d_{\mathcal V}(P(A),P(B))<\varepsilon/2.
\end{equation}
Combining \eqref{eq:prodsmall} and \eqref{eq:vaguesmall} yields
\[
d_{\sum}(\phi_*A,\phi_*B)
= d_{\mathrm{prod}}(\phi_*A,\phi_*B) + d_{\mathcal V}(P(\phi_*A),P(\phi_*B))
<\varepsilon.
\]
Thus $\phi_*$ is continuous.

Finally, the same argument applies to $\phi^{-1}\in\Aut(\N)$, so $(\phi^{-1})_*$ is continuous.
Since $(\phi^{-1})_*=(\phi_*)^{-1}$, it follows that $\phi_*$ is a homeomorphism.
\end{proof}

\subsection{Embedding $\pi_1(Conf_n(\C))$ into $\pi_1(\Conf)$}

\begin{definition}
Fix $n\ge 1$. For $x=(z_1,\dots,z_n)\in Conf_n(\C)$ put
\[
R(x):=1+\max_{1\le j\le n}|z_j|\in \R_{>0}.
\]
Define
\[
i_n:Conf_n(\C)\longrightarrow \Conf,\qquad
i_n(z_1,\dots,z_n):=(z_1,\dots,z_n,\ R(x),\ R(x)+1,\ R(x)+2,\dots).
\]
\end{definition}

\begin{lemma}
For each $n\ge 1$, the map $i_n$ is well-defined and continuous.
Moreover, $i_n(1,2,\dots,n)=\widetilde{\N}$.
\end{lemma}

\begin{proof}
Let $x=(z_1,\dots,z_n)\in Conf_n(\C)$ and set $R=R(x)$. For each $j\le n$ we have $|z_j|\le R-1<R$,
so $z_j\neq R+m$ for every integer $m\ge 0$. Thus there are no collisions between the first $n$
points and the appended tail. The tail points $(R+m)_{m\ge 0}$ are distinct and escape to $+\infty$
along the real axis, hence the resulting sequence is locally finite.

Continuity: the function $x\mapsto R(x)=1+\max_j|z_j|$ is continuous on $Conf_n(\C)$ (as a maximum of
finitely many continuous functions). If $x,x'\in Conf_n(\C)$ and $R=R(x),R'=R(x')$, then for each
tail coordinate we have $|(R+m)-(R'+m)|=|R-R'|$, so
\[
d_{\sum}\bigl(i_n(x),i_n(x')\bigr)\leq \sum_{j=1}^n|z_j-z'_j|+|R-R'|,
\]
which tends to $0$ as $x'\to x$. Finally, for $(1,2,\dots,n)$ we have $R=n+1$, so the tail is
$(n+1,n+2,\dots)$ and thus $i_n(1,2,\dots,n)=\widetilde{\N}$.
\end{proof}

\begin{definition}
Define
\[
p_n:\Conf\longrightarrow Conf_n(\C),\qquad p_n((a_j)_{j\ge 1}):=(a_1,\dots,a_n).
\]
\end{definition}

\begin{lemma}[Left inverse]
The map $p_n$ is continuous and satisfies $p_n\circ i_n=\mathrm{id}_{Conf_n(\C)}$.
\end{lemma}

\begin{proof}
Note that for $u=(u_1, ..., u_n), v=(v_1, ..., v_n)\in \C^n$,
$$d_n(u,v):=\sum_{j=1}^n 2^{-j}\min(|u_j-v_j|,1)$$
is a metric on $\C^n$ and it induces the usual product topology on \(\C^n\), hence also on
\(Conf_n(\C)\subset \C^n\).
Continuity of $p_n$ is immediate from the calculation
$$d_n(p_n(x), p_n(y))\leq d_{\mathrm{prod}}(x, y)\leq d_{\sum}(x, y)$$
The identity $p_n\circ i_n=\mathrm{id}$ holds because $i_n$ leaves the first $n$ coordinates unchanged.
\end{proof}

\begin{theorem}
For each $n\ge 1$, the induced homomorphism
\[
(i_n)_*:\pi_1\bigl(Conf_n(\C),(1,2,\dots,n)\bigr)\longrightarrow \pi_1(\Conf,\widetilde{\N})
\]
is injective. Equivalently, $\pi_1(Conf_n(\C))$ embeds as a subgroup of $\pi_1(\Conf)$.
\end{theorem}

\begin{proof}
Applying $\pi_1$ to $p_n\circ i_n=\mathrm{id}$ gives
\[
(p_n)_*\circ (i_n)_*=\mathrm{id}_{\pi_1(Conf_n(\C),(1,2,\dots,n))}.
\]
Hence $(i_n)_*$ has a left inverse and is therefore injective.
\end{proof}

\subsection{Uncountability of $\pi_1(\Conf)$}

\begin{lemma}\label{lem:uncountable_loop}
For each $k\ge 1$, let $c_k:=2k+\tfrac12$ and define a loop $\beta_k:[0,1]\to\Conf$ by
\[
\bigl(\beta_k(t)\bigr)_n=
\begin{cases}
c_k-\tfrac12 e^{2\pi i t}, & n=2k,\\[2pt]
c_k+\tfrac12 e^{2\pi i t}, & n=2k+1,\\[2pt]
n, & \text{otherwise.}
\end{cases}
\]
Then $\beta_k$ is a continuous loop based at $\widetilde{\N}$.
\end{lemma}

\begin{proof}
Set
\[
u(t):=c_k-\tfrac12 e^{2\pi it},\qquad v(t):=c_k+\tfrac12 e^{2\pi it}.
\]
The only moving coordinates are $2k$ and $2k+1$, and they remain distinct for all $t$ since
\[
|v(t)-u(t)|=\bigl|e^{2\pi it}\bigr|=1.
\]
All other coordinates are fixed at integers, so no collisions occur.
For each $t$ we have
\[
P(\beta_k(t))=\bigl(\N\setminus\{2k,2k+1\}\bigr)\ \cup\ \{u(t),v(t)\}.
\] which is locally finite.

Fix $t_0\in[0,1]$. For the product part,
\[
d_{\mathrm{prod}}(\beta_k(t),\beta_k(t_0))
=2^{-2k}\min\bigl(|u(t)-u(t_0)|,1\bigr)+2^{-(2k+1)}\min\bigl(|v(t)-v(t_0)|,1\bigr)\xrightarrow[t\to t_0]{}0,
\]
because $u$ and $v$ are continuous.

For the vague part, for each $j\ge1$ define
\[
S_j(t):=\sum_{a\in P(\beta_k(t))}\varphi_j(a).
\]
Since $\varphi_j$ has compact support and $P(\beta_k(t))$ is locally finite, the sum is finite.
Moreover, the only $t$-dependence of $S_j(t)$ is through $u(t)$ and $v(t)$, hence
\[
S_j(t)-S_j(t_0)=\bigl(\varphi_j(u(t))-\varphi_j(u(t_0))\bigr)+\bigl(\varphi_j(v(t))-\varphi_j(v(t_0))\bigr)\xrightarrow[t\to t_0]{}0.
\]
Then
\[
d_{\mathcal V}\bigl(P(\beta_k(t)),P(\beta_k(t_0))\bigr)
=\sum_{j=1}^{\infty}2^{-j}\frac{|S_j(t)-S_j(t_0)|}{1+|S_j(t)-S_j(t_0)|}\xrightarrow[t\to t_0]{}0.
\]
Therefore
\[
d_{\sum}\bigl(\beta_k(t),\beta_k(t_0)\bigr)
=d_{\mathrm{prod}}\bigl(\beta_k(t),\beta_k(t_0)\bigr)+d_{\mathcal V}\bigl(P(\beta_k(t)),P(\beta_k(t_0))\bigr)\xrightarrow[t\to t_0]{}0,
\]
so $\beta_k$ is continuous.
\end{proof}

\begin{definition}
For a binary sequence $\sigma=(\sigma_k)_{k\ge 1}\in\{0,1\}^{\N}$, define $\Gamma_\sigma:[0,1]\to\Conf$ by
\[
\bigl(\Gamma_\sigma(t)\bigr)_n=
\begin{cases}
c_k-\tfrac12 e^{2\pi i t}, & n=2k \text{ and } \sigma_k=1,\\[2pt]
c_k+\tfrac12 e^{2\pi i t}, & n=2k+1 \text{ and } \sigma_k=1,\\[2pt]
n, & \text{otherwise.}
\end{cases}
\]
\end{definition}

\begin{lemma}
For every $\sigma\in\{0,1\}^{\N}$, the map $\Gamma_\sigma$ is a continuous loop in $\Conf$ with
$\Gamma_\sigma(0)=\Gamma_\sigma(1)=\widetilde{\N}$.
\end{lemma}

\begin{proof}
The verification is the same as in Lemma \ref{lem:uncountable_loop}.
\end{proof}

\begin{lemma}
For each $k\ge 1$, the map
\[
F_k:\Conf\longrightarrow \C^*,\qquad F_k((a_j)_{j\geq 1}):=a_{2k}-a_{2k+1}
\]
is continuous and induces a homomorphism
\[
(F_k)_*:\pi_1(\Conf,\widetilde{\N})\longrightarrow \pi_1(\C^*,-1)\cong \mathbb{Z}.
\]
Moreover,
\[
(F_k)_*\bigl([\Gamma_\sigma]\bigr)=\sigma_k\ \ \text{for all }\sigma\in\{0,1\}^{\N}.
\]
\end{lemma}

\begin{proof}
Continuity of $F_k$ follows from continuity of coordinate projections in the metric $d_{\sum}$.
Since $a_{2k}\neq a_{2k+1}$ for every configuration in $\Conf$, the image lies in $\C^*$.

For the computation, if $\sigma_k=0$ then the coordinates $2k,2k+1$ are fixed at $2k,2k+1$, so
$F_k(\Gamma_\sigma(t))\equiv -1$ and the induced class $(F_k)_*([\Gamma_{\sigma}])$ is $0\in \Z$.

If $\sigma_k=1$, then
\[
F_k(\Gamma_\sigma(t))=(c_k-\tfrac12 e^{2\pi i t})-(c_k+\tfrac12 e^{2\pi i t})=-e^{2\pi i t},
\]
which winds once around $0$, hence represents $1\in\Z$ (for the standard orientation). Therefore
$(F_k)_*([\Gamma_\sigma])=\sigma_k$.
\end{proof}

\begin{theorem}
The fundamental group $\pi_1(\Conf,\widetilde{\N})$ is uncountable.
\end{theorem}

\begin{proof}
Define $\Phi:\{0,1\}^{\N}\to \pi_1(\Conf,\widetilde{\N})$ by $\Phi(\sigma)=[\Gamma_\sigma]$.
If $\sigma\neq\sigma'$, choose $k$ with $\sigma_k\neq\sigma'_k$. By Lemma~2.4,
\[
(F_k)_*(\Phi(\sigma))=\sigma_k\neq \sigma'_k=(F_k)_*(\Phi(\sigma')).
\]
Hence $\Phi(\sigma)\neq \Phi(\sigma')$, so $\Phi$ is injective. Since $\{0,1\}^{\N}$ has cardinality
$2^{\aleph_0}$, the group $\pi_1(\Conf,\widetilde{\N})$ is uncountable.
\end{proof}

\section{Path-connectedness of $\Conf$ and $\Clf$}

\subsection{A one-point isotopy supported in an open set}

\begin{lemma}[A bump homeomorphism moving a point a short distance]\label{lem:bumpmove}
Let $a\in\C$ and $r>0$. There exists a Lipschitz function $\eta:[0,\infty)\to[0,1]$ such that
$\eta\equiv1$ on $[0,r/3]$ and $\eta\equiv0$ on $[2r/3,\infty)$.
For any $p\in\C$ with $|p-a|\le r/6$, define
\[
h_{a,p}(z):=z+(p-a)\,\eta(|z-a|).
\]
Then $h_{a,p}$ is a homeomorphism of $\C$, satisfies $h_{a,p}(a)=p$, and is the identity on $\C\setminus B(a,2r/3)$.
Moreover, the map $(a,p,z)\mapsto h_{a,p}(z)$ is continuous on the set
\[
\{(a,p,z)\in\C^3:\ |p-a|\le r/6\}.
\]
In particular, if $a=a(x)$ and $p=p(x,t)$ depend continuously on parameters, then
$(x,t,z)\mapsto h_{a(x),p(x,t)}(z)$ is continuous.

\end{lemma}

\begin{proof}
Choose $\eta$ piecewise linear on $[0,r/3]$, $[r/3,2r/3]$, $[2r/3,\infty)$; then $\eta$ is Lipschitz with
$\mathrm{Lip}(\eta)\le 3/r$.

Fix $p$ with $|p-a|\le r/6$ and set $v:=p-a$. Define $g(z):=v\,\eta(|z-a|)$, so $h_{a,p}(z)=z+g(z)$.
Since $|g(z)-g(w)|\le |v|\,\mathrm{Lip}(\eta)\,|z-w|$, we have
\[
\mathrm{Lip}(g)\le |v|\,\mathrm{Lip}(\eta)\le (r/6)\cdot(3/r)=1/2.
\]

\smallskip
\noindent\emph{Injectivity and bi-Lipschitz bounds.}
For any $z,w\in\C$,
\[
|h_{a,p}(z)-h_{a,p}(w)|
=|(z-w)+(g(z)-g(w))|
\ge |z-w|-|g(z)-g(w)|
\ge (1-\mathrm{Lip}(g))|z-w|
\ge \tfrac12|z-w|.
\]
Thus $h_{a,p}$ is injective and the inverse (on its image) is $2$-Lipschitz.
Also
\[
|h_{a,p}(z)-h_{a,p}(w)|
\le |z-w|+|g(z)-g(w)|
\le (1+\mathrm{Lip}(g))|z-w|
\le \tfrac32|z-w|,
\]
so $h_{a,p}$ is globally Lipschitz.

\smallskip
\noindent\emph{Surjectivity (contraction mapping).}
Fix $y\in\C$. Consider the map $T_y:\C\to\C$ defined by
\[
T_y(z):=y-g(z).
\]
Then $T_y$ is a contraction with constant $\mathrm{Lip}(T_y)=\mathrm{Lip}(g)\le 1/2$.
By the Banach fixed point theorem, $T_y$ has a unique fixed point $z_y\in\C$, i.e.\ $z_y=y-g(z_y)$,
equivalently $h_{a,p}(z_y)=y$. Hence $h_{a,p}$ is surjective.

Therefore $h_{a,p}$ is bijective, continuous, and has a (globally) Lipschitz inverse, so it is a homeomorphism.
By construction $\eta\equiv 1$ near $0$, so $h_{a,p}(a)=a+v=p$, and $\eta\equiv 0$ on $[2r/3,\infty)$, so
$h_{a,p}=\id$ on $\C\setminus B(a,2r/3)$.
Continuity of $(p,z)\mapsto h_{a,p}(z)$ follows from continuity of $v=p-a$ and of $(z\mapsto \eta(|z-a|))$.
\end{proof}

\begin{lemma}[One-point isotopy supported in an open set]\label{lem:isotopy}
Let $U\subset\C$ be open and let $\gamma:[0,1]\to U$ be continuous. Then there exists a continuous family of
homeomorphisms $\Psi_t:\C\to\C$ $(t\in[0,1])$ such that:
\begin{enumerate}
\item $\Psi_0=\id$;
\item $\Psi_t(\gamma(0))=\gamma(t)$ for all $t$;
\item $\Psi_t=\id$ on $\C\setminus U$ for all $t$;
\item $(t,z)\mapsto \Psi_t(z)$ is continuous on $[0,1]\times\C$.
\end{enumerate}
Consequently, every $\Psi_t$ fixes every point of $\C\setminus U$.
\end{lemma}

\begin{proof}
Since $\gamma([0,1])$ is compact and contained in open $U$, there exists $\delta>0$ such that the closed $\delta$-neighborhood
of $\gamma([0,1])$ is contained in $U$.
By uniform continuity of $\gamma$, choose a partition $0=t_0<t_1<\cdots<t_N=1$ such that
$|\gamma(t)-\gamma(s)|<\delta/6$ whenever $t,s\in[t_{i-1},t_i]$.

Set $a_i:=\gamma(t_{i-1})$ and $r:=\delta$. For $t\in[t_{i-1},t_i]$ define
\[
\Phi_i(t):=h_{a_i,\gamma(t)}
\]
as in Lemma~\ref{lem:bumpmove}. Then $\Phi_i(t)$ is supported in $B(a_i,2r/3)\subset U$ and satisfies
$\Phi_i(t)(a_i)=\gamma(t)$.

Define $\Psi_t$ inductively: set $\Psi_{t_0}:=\id$, and for $t\in[t_{i-1},t_i]$ put
\[
\Psi_t:=\Phi_i(t)\circ \Psi_{t_{i-1}}.
\]
Since $\Phi_i(t_{i-1})=\id$, the definition matches at the endpoints and gives continuity of $t\mapsto\Psi_t$.
Each $\Psi_t$ is the identity outside $U$ because each $\Phi_i(t)$ is.
Finally, by induction on $i$ we have $\Psi_{t_{i-1}}(\gamma(0))=\gamma(t_{i-1})=a_i$, hence for $t\in[t_{i-1},t_i]$,
\[
\Psi_t(\gamma(0))=\Phi_i(t)(a_i)=\gamma(t).
\]
\end{proof}

\subsection{Single-step stabilization with escaping support}

\begin{lemma}\label{lem:detour}
Let $S\subset\C$ be locally finite and let $\sigma:[0,1]\to\C$ be continuous with compact image.
Then $\sigma([0,1])\cap S$ is finite. Moreover, if $\sigma(0),\sigma(1)\notin S$, there exists a continuous path
$\tilde\sigma:[0,1]\to \C\setminus S$ with $\tilde\sigma(0)=\sigma(0)$ and $\tilde\sigma(1)=\sigma(1)$
whose image is contained in an arbitrarily small neighborhood of $\sigma([0,1])$.
\end{lemma}
\begin{proof}
Since $\sigma([0,1])$ is compact and $S$ is locally finite, $\sigma([0,1])\cap S$ is finite.
Let $F=\sigma([0,1])\cap S$.
Choose pairwise disjoint closed disks $\cl{B(p,\rho_p)}$ around each $p\in F$ such that
$\cl{B(p,\rho_p)}\cap (S\setminus\{p\})=\emptyset$ (possible since $S$ is closed and $F$ is finite).
Modify $\sigma$ inside each disk by replacing the (finitely many) subarcs that hit $p$ with a small detour in
$\cl{B(p,\rho_p)}\setminus\{p\}$ connecting entry/exit points. This yields $\tilde\sigma$ as required.
\end{proof}

\begin{lemma}[Single-step stabilization with clearing of the target]\label{lem:singlestep}
Let $n\ge1$ and let $y=(y_j)_{j\geq 1}\in \Conf$ satisfy $y_j=j$ for $1\le j\le n-1$.
Put
\[
S_n:=\{y_j:\ j\neq n\}.
\]
Define $m_n:=\min\{|y_n|,\,n\}$ and
\[
R_n:=\frac{m_n}{4}-2.
\]
Then there exist:
\begin{itemize}
\item an index set $J_n\subset\N$ with $n\in J_n$ and $|J_n|\le 2$;
\item an open set $U_n\subset\{z:\ |z|>R_n\}$ with $U_n\cap\{y_j:\ j\notin J_n\}=\emptyset$;
\item a continuous path $G^{(n)}:[0,1]\to (\Conf, d_{\mathrm{prod}})$;
\end{itemize}
such that, writing $g^{(n)}:=G^{(n)}(1)$:
\begin{enumerate}[label=(\roman*)]
\item $G^{(n)}(0)=y$ and $g^{(n)}_j=j$ for $1\le j\le n$;
\item for every $t\in[0,1]$ and every $j\notin J_n$, one has $G^{(n)}_j(t)=y_j$;
\item for every $t\in[0,1]$ and every $j\in J_n$, one has $|G^{(n)}_j(t)|\ge R_n+1$.
\end{enumerate}
Moreover, if $J_n=\{n,k\}$, then $k>n$ is the unique index with $y_k=n$, and the construction can be arranged so that
$|g^{(n)}_k|\ge k$.
\end{lemma}

\begin{proof}
\textbf{Step 1 (endpoint margin).}
Since $R_n+2=m_n/4\le m_n\le |y_n|$ and $R_n+2\le m_n\le n$, both points $y_n$ and $n$ lie in $\{|z|\ge R_n+2\}$.

\medskip
\noindent\textbf{Step 2 (Case 1: $n\notin S_n$).}
Then the endpoint $n$ is not occupied by any other coordinate.
Set $J_n=\{n\}$.
Choose a polygonal path $\sigma:[0,1]\to \{|z|\ge R_n+2\}$ with $\sigma(0)=y_n$ and $\sigma(1)=n$.
By Lemma~\ref{lem:detour}, modify $\sigma$ to a path
\[
\gamma:[0,1]\to \{|z|\ge R_n+1\}\setminus S_n
\]
(from $y_n$ to $n$); we may require $\gamma$ to stay within distance $1/2$ of $\sigma$, hence in $\{|z|\ge R_n+1\}$.
Since $\gamma([0,1])$ is compact and disjoint from $S_n$,
the distance $\eta:=\dist(\gamma([0,1]),S_n)>0$.
Let $\rho:=\tfrac12\min\{1,\eta\}$ and define the tubular neighborhood
\[
U_n:=\{z\in\C:\dist(z,\gamma([0,1]))<\rho\}.
\]
Then $U_n\subset\{|z|>R_n\}$ (since $\gamma$ lies in $\{|z|\ge R_n+1\}$ and $\rho\le 1/2$), and $U_n\cap S_n=\emptyset$.
Apply Lemma~\ref{lem:isotopy} to $U_n$ and $\gamma$ to obtain $\Psi_t$ supported in $U_n$ with $\Psi_t(y_n)=\gamma(t)$.
Define $G^{(n)}(t):=(\Psi_t(y_j))_{j\ge1}$. Then only the $n$-th coordinate moves, and it stays in $\{|z|\ge R_n+1\}$.

\medskip
\noindent\textbf{Step 3 (Case 2: $n\in S_n$).}
Then there is a unique index $k\neq n$ with $y_k=n$. Since $y_j=j$ for $j\le n-1$, necessarily $k>n$.
Set $J_n:=\{n,k\}$.
We first clear the target $n$ by moving the $k$-th coordinate away, keeping all other coordinates fixed.

\smallskip
\noindent\emph{Step 3a (choose a far new location for the blocker).}
Let $\rho_{n,k}:=\max\{R_n+3,\,k\}$ and let $C_{\rho_{n,k}}:=\{z:\ |z|=\rho_{n,k}\}$.
Since $P(y)$ is locally finite, $C_{\rho_{n,k}}\cap P(y)$ is finite, hence we can choose
\[
q\in C_{\rho_{n,k}}\setminus P(y).
\]
Note that $|q|=\rho_{n,k}\ge k$.

\smallskip
\noindent\emph{Step 3b (move the blocker from $n$ to $q$).}
Let $S_k:=\{y_j:\ j\neq k\}$, a locally finite (hence closed) set that does not contain $n$.
Choose a polygonal path $\sigma:[0,1]\to\{|z|\ge R_n+2\}$ from $\sigma(0)=n$ to $\sigma(1)=q$.
By Lemma~\ref{lem:detour}, modify $\sigma$ to a path
\[
\alpha:[0,1]\to \{|z|\ge R_n+1\}\setminus S_k
\]
with $\alpha(0)=n$ and $\alpha(1)=q$.
Let $\eta:=\dist(\alpha([0,1]),S_k)>0$ and put $\rho:=\tfrac12\min\{1,\eta\}$.
Set
\[
V:=\{z:\dist(z,\alpha([0,1]))<\rho\}\subset\{|z|>R_n\}.
\]
Then $V\cap S_k=\emptyset$. Apply Lemma~\ref{lem:isotopy} to obtain homeomorphisms $\Theta_t$ supported in $V$
with $\Theta_t(n)=\alpha(t)$ and fixing $\C\setminus V$ pointwise.
Define an intermediate configuration $\tilde y$ by $\tilde y_j:=\Theta_1(y_j)$.
Then $\tilde y_k=q$ and $\tilde y_j=y_j$ for $j\neq k$ (in particular $\tilde y_n=y_n$).

\smallskip
\noindent\emph{Step 3c (move coordinate $n$ from $y_n$ to $n$, now unoccupied).}
Let $\tilde S_n:=\{\tilde y_j:\ j\neq n\}$. Since $\tilde y_k=q\neq n$, we have $n\notin \tilde S_n$.
Choose a polygonal path $\sigma':[0,1]\to\{|z|\ge R_n+2\}$ from $\sigma'(0)=y_n$ to $\sigma'(1)=n$ and,
by Lemma~\ref{lem:detour}, modify it to a path
\[
\beta:[0,1]\to \{|z|\ge R_n+1\}\setminus \tilde S_n.
\]
Let $\eta':=\dist(\beta([0,1]),\tilde S_n)>0$ and set $\rho':=\tfrac12\min\{1,\eta'\}$.
Let
\[
W:=\{z:\dist(z,\beta([0,1]))<\rho'\}\subset\{|z|>R_n\}.
\]
Then $W\cap \tilde S_n=\emptyset$. Apply Lemma~\ref{lem:isotopy} to obtain homeomorphisms $\Psi_t$ supported in $W$
with $\Psi_t(y_n)=\beta(t)$ and fixing $\C\setminus W$ pointwise; in particular $\Psi_t$ fixes $q$.

\smallskip
\noindent\emph{Step 3d (concatenate the two moves).}
Define $G^{(n)}:[0,1]\to \Conf$ by
\[
G^{(n)}(t):=
\begin{cases}
\bigl(\Theta_{2t}(y_j)\bigr)_{j\ge1}, & 0\le t\le \tfrac12,\\[3pt]
\bigl(\Psi_{2t-1}(\tilde y_j)\bigr)_{j\ge1}, & \tfrac12\le t\le 1.
\end{cases}
\]
Then $G^{(n)}(0)=y$, and $G^{(n)}(1)$ satisfies $g^{(n)}_n=n$. Only indices in $J_n=\{n,k\}$ move.
Moreover, the moving coordinates stay in $\{|z|\ge R_n+1\}$, and the final blocker position $g^{(n)}_k=q$ satisfies $|q|\ge k$.

\smallskip
Finally, set $U_n:=U_n$ in Case~1, and $U_n:=V\cup W$ in Case~2. In both cases $U_n\subset\{|z|>R_n\}$ and
$U_n\cap\{y_j:\ j\notin J_n\}=\emptyset$.

\medskip
\noindent\textbf{Step 4 (continuity and membership in $\Conf$).}
Each $G^{(n)}(t)$ is obtained from $y$ by a homeomorphism of $\C$ (or concatenation of two),
hence has pairwise distinct coordinates and locally finite underlying set; thus $G^{(n)}(t)\in \Conf$.
Only finitely many coordinates move (at most two), and each moving coordinate depends continuously on $t$,
so $G^{(n)}$ is continuous in $d_{\mathrm{prod}}$.
\end{proof}

\subsection{Iterated stabilization and a global path}

\begin{lemma}[Iterated stabilization]\label{lem:iterated}
Fix $x\in \Conf$. There exist configurations $g^{(n)}\in \Conf$ ($n\ge0$), paths $G^{(n)}:[0,1]\to \Conf$ ($n\ge1$),
and radii $R_n\in\R$ such that:
\begin{enumerate}[label=(\roman*)]
\item $g^{(0)}=x$;
\item $g^{(n)}_j=j$ for $1\le j\le n$ (hence those coordinates never change again);
\item $G^{(n)}(0)=g^{(n-1)}$ and $G^{(n)}(1)=g^{(n)}$;
\item for each $n$, there is a set $J_n\subset\N$ with $|J_n|\le 2$ such that for all $t$ and all $j\notin J_n$,
\ $G^{(n)}_j(t)=g^{(n-1)}_j$;
\item the paths $G^{(n)}$ are induced by isotopies supported in $\{|z|>R_n\}$ and the moving coordinates satisfy
$|G^{(n)}_j(t)|\ge R_n+1$ for $j\in J_n$; moreover $R_n\to\infty$ as $n\to\infty$.
\end{enumerate}
\end{lemma}
\begin{proof}
Induct on $n$, applying Lemma~\ref{lem:singlestep} to $y=g^{(n-1)}$.
This produces $G^{(n)}$, $J_n$, and $R_n$ with properties (ii)--(iv).

It remains to prove $R_n\to\infty$.
In Lemma~\ref{lem:singlestep}, $R_n$ is defined as $R_n=\frac{1}{4}\min\{|g^{(n-1)}_n|,n\}-2$.
We claim that for each $n$,
\[
|g^{(n-1)}_n|\ \ge\ \min\{|x_n|,\,n\}.
\]
Indeed, either coordinate $n$ was never moved before stage $n$, in which case $g^{(n-1)}_n=x_n$,
or it was moved earlier as a blocker in some stage $m<n$; in that case Lemma~\ref{lem:singlestep} guarantees that
after that move its modulus is at least its index, hence at least $n$.
Thus $\min\{|g^{(n-1)}_n|,n\}\ge \min\{|x_n|,n\}$, and therefore
\[
R_n \ \ge\ \frac{1}{4}\min\{|x_n|,n\}-2 \ \xrightarrow[n\to\infty]{}\ \infty
\]
\end{proof}

\begin{definition}[Stage intervals and global path]\label{def:global}
Let
\[
I_1=\Bigl[0,\frac12\Bigr],\qquad
I_n=\bigl(1-2^{-(n-1)},\,1-2^{-n}\bigr]\quad(n\ge2),
\]
so that $[0,1)=\bigsqcup_{n\ge1}I_n$.
Let $\theta_n:I_n\to[0,1]$ be the affine increasing bijection.

Given $(g^{(n)})$ and $(G^{(n)})$ from Lemma~\ref{lem:iterated}, define $H:[0,1]\to \Conf$ by
\[
H(t)=
\begin{cases}
G^{(n)}(\theta_n(t)), & t\in I_n,\\
\widetilde{\N}, & t=1.
\end{cases}
\]
\end{definition}

\begin{lemma}[Continuity of $H$ in $d_{\mathrm{prod}}$]\label{lem:Hprod}
The map $H:[0,1]\to(E,d_{\mathrm{prod}})$ is continuous.
\end{lemma}
\begin{proof}
On each $I_n$, $H=G^{(n)}\circ\theta_n$ is continuous.
At the junction time $t_n:=\sup I_n=1-2^{-n}$ we have
\[
\lim_{t\nearrow t_n}H(t)=G^{(n)}(1)=g^{(n)}=G^{(n+1)}(0)=\lim_{t\searrow t_n}H(t),
\]
so $H$ is continuous on $[0,1)$.

For continuity at $1$, fix $\varepsilon>0$ and choose $M$ so that $\sum_{j>M}2^{-j}<\varepsilon$.
Once stage $j$ is completed, coordinate $j$ equals $j$ forever (Lemma~\ref{lem:iterated}(ii)),
so there exists $T<1$ such that $H_j(t)=j=H_j(1)$ for all $1\le j\le M$ and all $t\in[T,1]$.
Hence for $t\in[T,1]$,
\[
d_{\mathrm{prod}}(H(t),H(1))
=\sum_{j>M}2^{-j}\min(|H_j(t)-j|,1)\le \sum_{j>M}2^{-j}<\varepsilon.
\]
\end{proof}

For $\varphi\in C_c(\C)$ define the test sum
\[
S_\varphi(t):=\sum_{a\in P(H(t))}\varphi(a),
\]
which is a finite sum.

\begin{lemma}[Continuity of test sums]\label{lem:testsum}
For each fixed $\varphi\in C_c(\C)$, the map $t\mapsto S_\varphi(t)$ is continuous on $[0,1]$.
\end{lemma}
\begin{proof}
Let $K=\supp(\varphi)$, compact, and choose $R>0$ with $K\subset \cl{B(0,R)}$.

\smallskip
\noindent\textbf{Continuity at $t_0<1$.}
Fix $t_0<1$ and choose a compact interval $J\subset[0,1)$ containing $t_0$.
Then $J$ meets only finitely many stage intervals $I_n$, hence only finitely many stages occur on $J$.
Let $F\subset\N$ be the union of the corresponding moving index sets $J_n$; then $F$ is finite, and only indices in $F$
can move on $J$.

For indices $j\notin F$, the point $H_j(t)$ is constant on $J$.
Since $P(H(t_0))$ is locally finite, $P(H(t_0))\cap K$ is finite.
Therefore there exists a finite index set $I\subset\N$ such that for all $t\in J$,
\[
P(H(t))\cap K \subset \{H_i(t): i\in I\}.
\]
(Indeed, take $I:=F\cup\{j\notin F : H_j(t_0)\in K\}$, which is finite because $P(H(t_0))\cap K$ is finite.)
Hence on $J$ we have $S_\varphi(t)=\sum_{i\in I}\varphi(H_i(t))$, a finite sum of continuous functions, so
$S_\varphi$ is continuous at $t_0$.

\smallskip
\noindent\textbf{Continuity at $t_0=1$.}
By Lemma~\ref{lem:iterated}(v), the radii $R_n\to\infty$, so choose $N_1$ such that $R_n>R$ for all $n\ge N_1$.
For such $n$, every moving coordinate during stage $n$ stays in $\{|z|\ge R_n+1\}\subset \C\setminus \cl{B(0,R)}$,
hence cannot enter $K$.

Let $I_0:=\{j\in\N : x_j\in K\}$, finite because $P(x)$ is locally finite.
Let $I_1:=\bigcup_{n=1}^{N_1-1} J_n$, finite (each $J_n$ has size $\le2$).
Set $I:=I_0\cup I_1$, finite, and let $N:=\max I$.

\emph{Claim:} for every $t\in[0,1]$, $P(H(t))\cap K \subset \{H_i(t): i\in I\}$.
Indeed, if $H_j(t)\in K$ for some $j$, then either $j$ never moved up to time $t$ (so $H_j(t)=x_j\in K$ and $j\in I_0$),
or $j$ moved during some stage $n$ with $n<N_1$ (hence $j\in J_n\subset I_1$),
because stages $n\ge N_1$ never place moving coordinates into $K$. In other words,
a point can lie in $K$ only if it either started in $K$ and never moved $(I_0)$, or it was one of the finitely many indices moved during the finitely many “early” stages where motion might still happen near $K(I_1)$; after stage
$N_1$, all motion is forced far outside $K$.

Now take $t\ge t_N:=1-2^{-N}$. Then all stages $1,\dots,N$ have completed, so $H_i(t)=i$ for every $i\le N$.
In particular $H_i(t)=i$ for all $i\in I$.
If $j>N$, then $j\notin I$, hence $H_j(t)\notin K$ by the claim above.
Therefore for all $t\in[t_N,1]$,
\[
P(H(t))\cap K = \{1,2,3,\dots\}\cap K = P(\widetilde{\N})\cap K,
\]
so $S_\varphi(t)=S_\varphi(1)$ for all $t\in[t_N,1]$, proving continuity at $1$.
\end{proof}

\begin{lemma}[Continuity of $P\circ H$ in $d_{\mathcal V}$]\label{lem:Hvague}
The map $P\circ H:[0,1]\to(C^{lf}_\infty(\C),d_{\mathcal V})$ is continuous.
\end{lemma}
\begin{proof}
Fix $t_0\in[0,1]$.
For each $m$, Lemma~\ref{lem:testsum} gives continuity of $t\mapsto S_{\varphi_m}(t)$ at $t_0$.
Define
\[
f_m(t):=\frac1{2^m}\frac{|S_{\varphi_m}(t)-S_{\varphi_m}(t_0)|}{1+|S_{\varphi_m}(t)-S_{\varphi_m}(t_0)|}.
\]
Then $f_m$ is continuous at $t_0$ and $0\le f_m(t)\le 2^{-m}$.
By the Weierstrass $M$-test, $\sum_m f_m(t)$ converges uniformly in a neighborhood of $t_0$.
Hence
\[
t\longmapsto d_{\mathcal V}(P(H(t)),P(H(t_0)))=\sum_{m\ge1} f_m(t)
\]
is continuous at $t_0$.
\end{proof}

\begin{theorem}\label{thm:pathconn}
The metric space $(E,d_{\sum})$ is path-connected.
\end{theorem}

\begin{proof}
Construct $H$ as in Definition~\ref{def:global}. By Lemma~\ref{lem:Hprod}, $H$ is continuous into $(E,d_{\mathrm{prod}})$.
By Lemma~\ref{lem:Hvague}, $P\circ H$ is continuous into $(C^{lf}_\infty(\C),d_{\mathcal V})$.
Since $d_{\sum}(u,v)=d_{\mathrm{prod}}(u,v)+d_{\mathcal V}(P(u),P(v))$, it follows that $H$ is continuous into $(E,d_{\sum})$.
\end{proof}

Since $P: \Conf \rightarrow \Clf$ is continuous and surjective, we have the following result.

\begin{corollary}
The space $(\Clf, d_{\mathcal V})$ is path-connected.
\end{corollary}

\section{Asphericity of $\Clf$}

%%%%%%%%%%%%%%%%%%%%%%%%%%%%%%%%%%%%%%%%%%%%%%%%%%%%%%%%%%%%%%%%%%%%%%%%%%%%%%
%  INPUTS assumed available (as in your framework)
%%%%%%%%%%%%%%%%%%%%%%%%%%%%%%%%%%%%%%%%%%%%%%%%%%%%%%%%%%%%%%%%%%%%%%%%%%%%%%
% Lemma~\ref{lem:ub}: Uniform boundedness on compacta (vague topology).
% Lemma~\ref{lem:param-push}: Parametrized push-forward is continuous under fixed compact support.
% Lemma~\ref{lem:restrict}: Restriction lemma: if F(x) avoids ∂B for all x, then counts in B are locally constant
%                          and restriction to int(B) is continuous into finite configuration spaces.
% Theorem~\ref{thm:EK}: Edwards--Kirby parametric relative isotopy extension with support control.
% Theorem~\ref{thm:FA} + Lemma~\ref{lem:sc-to-kpi1}: Conf_n(disk/annulus) is K(pi,1), hence maps from simply connected CW are null-homotopic.

\subsection{Relative isotopy extension}
\begin{lemma}\label{lem:ub}
Let $K$ be a compact topological space and let
\[
F:K\longrightarrow C^{lf}_\infty(\C)
\]
be continuous. Then for every compact set $B\subset\C$
there exists an integer $N(B)\ge 0$ such that
\[
\#\bigl(F(x)\cap B\bigr)\le N(B)\qquad\text{for all }x\in K.
\]
Equivalently, the family $\{F(x):x\in K\}$ is uniformly locally finite on every compact subset of $\C$.
\end{lemma}

\begin{proof}
Choose $\psi\in C_c(\C)$ such that
\[
0\le \psi\le 1,\qquad \psi\equiv 1 \ \text{on }B.
\]
Define
\[
g:K\to \R,\qquad g(x):=\langle F(x),\psi\rangle=\sum_{a\in F(x)}\psi(a).
\]
This sum is finite for each $x$ because $F(x)$ is locally finite and $\supp(\psi)$ is compact.
By definition of the vague topology, $g$ is continuous. Since $K$ is compact, $g$ is bounded:
there exists $M\in \N$ such that $g(x)\le M$ for all $x\in K$.

Now note that for each $x\in K$,
\[
g(x)=\sum_{a\in F(x)}\psi(a)\ \ge\ \sum_{a\in F(x)\cap B}\psi(a)
\ =\ \sum_{a\in F(x)\cap B}1
\ =\ \#(F(x)\cap B),
\]
because $\psi\equiv 1$ on $B$ and $\psi\ge 0$ everywhere. Hence
\[
\#(F(x)\cap B)\le g(x)\le M\qquad\mbox{ for all } x\in K.
\]
\end{proof}

\begin{lemma}\label{lem:param-push}
Let $K$ be a compact topological space and let $F:K\to C^{lf}_\infty(\C)$ be continuous. Let
\[
\Phi:K\times[0,1]\times \C \longrightarrow \C
\]
be continuous, and write $\Phi_{x,t}:=\Phi(x,t,\cdot)$.
Assume:
\begin{enumerate}
\item for every $(x,t)\in K\times[0,1]$, the map $\Phi_{x,t}:\C\to\C$ is a homeomorphism;
\item there exists a fixed compact set $S\subset\C$ such that
\[
\Phi_{x,t}(z)=z\qquad \text{for all }(x,t)\in K\times[0,1]\text{ and all }z\in \C\setminus S.
\]
\end{enumerate}
Define
\[
H:K\times[0,1]\to C^{lf}_\infty(\C),\qquad H(x,t):=(\Phi_{x,t})_*\bigl(F(x)\bigr)
=\{\Phi_{x,t}(a):a\in F(x)\}.
\]
Then $H$ is continuous.
\end{lemma}

\begin{proof}
To prove that $H:K\times[0,1]\to C^{lf}_\infty(\C)$ is continuous in the vague topology,
it suffices to show that for every $\varphi\in C_c(\C)$ the map
\[
g(x,t):=\langle H(x,t),\varphi\rangle
=\sum_{a\in F(x)}\varphi\bigl(\Phi_{x,t}(a)\bigr)
\]
is continuous on $K\times[0,1]$.

Fix $\varphi\in C_c(\C)$ and set $C_0:=\supp(\varphi)$, which is compact.

\medskip
\noindent\textbf{Step 1}
Define
\[
C_1:=\bigcup_{(x,t)\in K\times[0,1]}\Phi_{x,t}^{-1}(C_0)\subset \C.
\]
We claim $C_1$ is compact.

Indeed, by the support assumption there exists a fixed compact $S\subset\C$ such that
$\Phi_{x,t}(z)=z$ for all $z\notin S$. Hence, if $z\in \Phi_{x,t}^{-1}(C_0)$ and $z\notin S$,
then $z=\Phi_{x,t}(z)\in C_0$. Therefore
$
C_1\subset S\cup C_0,
$
and $S\cup C_0$ is compact.

Now consider the set
\[
E:=\{(x,t,z)\in K\times[0,1]\times (S\cup C_0)\;:\;\Phi_{x,t}(z)\in C_0\}.
\]
Since $(x,t,z)\mapsto \Phi_{x,t}(z)$ is continuous, $E$ is closed in the compact space
$K\times[0,1]\times(S\cup C_0)$, hence $E$ is compact. Let $\pi_\C$ be the projection to the $z$-coordinate.
Then
$
C_1=\pi_\C(E),
$
so $C_1$ is compact as the continuous image of a compact set.

\medskip
\noindent\textbf{Step 2}
Since $C_1$ is compact and $F:K\to C^{lf}_\infty(\C)$ is continuous with $K$ compact,
Lemma~\ref{lem:ub} gives an integer $N\ge0$ such that
\[
\#\bigl(F(x)\cap C_1\bigr)\le N\qquad\mbox{ for all } x\in K.
\]
Moreover, $\varphi(\Phi_{x,t}(a))\neq 0$ implies $\Phi_{x,t}(a)\in C_0$, i.e.\ $a\in C_1$.
Hence we may rewrite
\[
g(x,t)=\sum_{a\in F(x)\cap C_1}\varphi\bigl(\Phi_{x,t}(a)\bigr).
\]

If $N=0$ then $F(x)\cap C_1=\varnothing$ for all $x$, so $g\equiv 0$ and we are done.
Assume henceforth $N\ge1$.

\medskip
\noindent\textbf{Step 3}
Fix $(x_0,t_0)\in K\times[0,1]$. Define, for $z\in C_1$,
\[
f_{x,t}(z):=\varphi\bigl(\Phi_{x,t}(z)\bigr),\qquad
f_0(z):=f_{x_0,t_0}(z)=\varphi\bigl(\Phi_{x_0,t_0}(z)\bigr).
\]
Then $(x,t,z)\mapsto |f_{x,t}(z)-f_0(z)|$ is continuous on the compact space
$K\times[0,1]\times C_1$. Therefore the function
\[
m(x,t):=\sup_{z\in C_1}|f_{x,t}(z)-f_0(z)|
\]
is continuous on $K\times[0,1]$.
In particular, $m(x_0,t_0)=0$, so there exists a neighborhood $U$ of $(x_0,t_0)$ such that
\[
m(x,t)<\frac{\varepsilon}{2N}\qquad\mbox{ for all } (x,t)\in U.
\tag{$\ast$}
\]

Next, note that $f_0\in C_c(\C)$ because $\Phi_{x_0,t_0}$ is a homeomorphism and
\[
\supp(f_0)\subset \Phi_{x_0,t_0}^{-1}(C_0)\subset C_1
\]
is compact. Hence the map $x\mapsto \langle F(x),f_0\rangle$ is continuous.
Thus there exists a neighborhood $U'\subset K$ of $x_0$ such that
\[
\bigl|\langle F(x),f_0\rangle-\langle F(x_0),f_0\rangle\bigr|<\frac{\varepsilon}{2}
\qquad\mbox{ for all } x\in U'.
\tag{$\ast\ast$}
\]

Now let $(x,t)\in (U'\times[0,1])\cap U$. We estimate:
\[
\begin{aligned}
|g(x,t)-g(x_0,t_0)|
&=\left|\sum_{a\in F(x)\cap C_1} f_{x,t}(a) - \sum_{b\in F(x_0)\cap C_1} f_0(b)\right|\\
&\le
\left|\sum_{a\in F(x)\cap C_1}\bigl(f_{x,t}(a)-f_0(a)\bigr)\right|
+
\left|\sum_{a\in F(x)\cap C_1} f_0(a) - \sum_{b\in F(x_0)\cap C_1} f_0(b)\right|\\
&\le \#(F(x)\cap C_1)\cdot m(x,t)
+
\bigl|\langle F(x),f_0\rangle-\langle F(x_0),f_0\rangle\bigr|\\
&\le N\cdot \frac{\varepsilon}{2N} + \frac{\varepsilon}{2}
=\varepsilon,
\end{aligned}
\]
using $(\ast)$ and $(\ast\ast)$. Therefore $g$ is continuous at $(x_0,t_0)$, and hence continuous on
$K\times[0,1]$.

\medskip
Since this holds for every $\varphi\in C_c(\C)$, the map $H$ is continuous in the vague topology.
\end{proof}

We first record a basic openness fact that will be used repeatedly.

\begin{lemma}\label{lem:openness}
Let $K$ be a compact topological space and $F: K \rightarrow \Clf$ be continuous.
If $Q\subset \C$ is compact and $F(x_0)\cap Q=\varnothing$,
then there exists an open neighborhood $U$ of $x_0$ in $K$ such that
$F(x)\cap Q=\varnothing$ for all $x\in U$.
\end{lemma}

\begin{proof}
Choose $\psi\in C_c(\C)$ with $\psi\ge 1$ on $Q$ and $\supp\psi\cap F(x_0)=\varnothing$.
Then $\langle F(x_0),\psi\rangle=0$. By vague continuity, $x\mapsto \langle F(x),\psi\rangle$
is continuous, hence there is a neighborhood $U$ of $x_0$ with
$\langle F(x),\psi\rangle<1/2$ for all $x\in U$.
If $F(x)\cap Q\neq\varnothing$, then $\langle F(x),\psi\rangle\ge 1$, contradiction.
\end{proof}

\begin{lemma}\label{lem:restrict}
Let $K$ be a compact topological space and let
\[
F:K\longrightarrow C^{lf}_\infty(\C)
\]
be continuous. Let $B\subset\C$ be a compact set such that
\[
F(x)\cap \partial B=\varnothing \qquad \text{for all }x\in K.
\]
Then:

\begin{enumerate}
\item The function
\[
x\longmapsto \#\bigl(F(x)\cap B\bigr)
\]
is locally constant on $K$.

\item For each connected component $K_0\subset K$, let
\[
n:=\#\bigl(F(x)\cap B\bigr)\quad(\text{ where }x\in K_0).
\]
Then the map
\[
r_B:K_0\longrightarrow C_n(\operatorname{int}B),\qquad
r_B(x):=F(x)\cap \operatorname{int}B
\]
is well-defined and continuous, where $C_n(\operatorname{int}B)$ denotes the n-th unordered configuration space of $\operatorname{int}B$.
\end{enumerate}
\end{lemma}

\begin{proof}

\noindent\textbf{(1) Local constancy of $x\mapsto \#(F(x)\cap B)$.}
Fix $x_0\in K$. Since $F(x_0)$ is locally finite, $F(x_0)\cap \partial B$ is finite; by assumption it is empty.
Because $\partial B$ is compact and $F(x_0)$ has no accumulation points in $\partial B$,
the distance
\[
\delta:=\dist(\partial B,F(x_0))>0.
\]
Let
\[
N:=\{z\in\C:\dist(z,\partial B)\le \delta/3\},
\]
a compact collar of $\partial B$. Then $F(x_0)\cap N=\varnothing$, hence by Lemma \ref{lem:openness},
there is a neighborhood $U$ of $x_0$ such that
\[
F(x)\cap N=\varnothing\qquad\mbox{ for all } x\in U.
\tag{$\ast$}
\]
Define the compact set
\[
B^-:=\{z\in B:\dist(z,\partial B)\ge \delta/3\}\subset \operatorname{int}B.
\]
For $x\in U$, condition ($\ast$) implies that $F(x)$ has no points in the strip
$B\setminus B^-\subset N$, hence
\[
F(x)\cap B=F(x)\cap B^- \qquad (x\in U).
\tag{$\ast\ast$}
\]
Choose $\eta\in C_c(\C)$ with $0\le \eta\le 1$, $\eta\equiv 1$ on $B^-$, and
$\supp(\eta)\subset B^-\cup N$ (possible because $B^-\subset \operatorname{int}B$ and $N$ is compact).
Then for $x\in U$, using ($\ast$) and ($\ast\ast$),
\[
\langle F(x),\eta\rangle
=\sum_{a\in F(x)}\eta(a)
=\sum_{a\in F(x)\cap B^-} 1
=\#(F(x)\cap B^-)
=\#(F(x)\cap B).
\]
Thus $x\mapsto \#(F(x)\cap B)$ agrees on $U$ with the continuous map
$x\mapsto \langle F(x),\eta\rangle$ and is therefore continuous on $U$.
But it takes values in $\N$, hence is locally constant on $U$.
Since $x_0$ was arbitrary, $x\mapsto \#(F(x)\cap B)$ is locally constant on $K$.

\medskip
\noindent\textbf{(2) Continuity of the restriction map $r_B$ on each component.}
Let $K_0\subset K$ be a connected component. By (1) the number
\[
n:=\#(F(x)\cap B)
\]
is constant for $x\in K_0$. Since $F(x)\cap \partial B=\varnothing$, we have
$F(x)\cap B\subset \operatorname{int}B$, so $r_B(x)=F(x)\cap \operatorname{int}B$ is well-defined
as an element of $C_n(\operatorname{int}B)$.

Fix $x_0\in K_0$.
If $n=0$, then $r_B$ is the constant map to the empty configuration and is continuous.
Assume $n\ge 1$.

As above, choose $\delta>0$ and a neighborhood $U\subset K_0$ of $x_0$ so that ($\ast$) holds and hence
$F(x)\cap B=F(x)\cap B^-$ for all $x\in U$. Set
\[
A_0:=F(x_0)\cap B^-=\{a_1,\dots,a_n\}\subset \operatorname{int}B.
\]
Choose pairwise disjoint open disks $D_1,\dots,D_n$ in $\C$ such that:
\[
a_i\in D_i,\qquad \overline{D_i}\subset \operatorname{int}B,
\qquad \overline{D_i}\cap \overline{D_j}=\varnothing\ (i\neq j),
\]
and also such that $F(x_0)\cap \overline{D_i}=\{a_i\}$ for each $i$.
Let
\[
C_{\mathrm{out}}:=B^-\setminus \bigcup_{i=1}^n D_i,
\]
which is compact and satisfies $F(x_0)\cap C_{\mathrm{out}}=\varnothing$.
By Claim~0, after shrinking $U$ we may assume
\[
F(x)\cap C_{\mathrm{out}}=\varnothing\qquad\mbox{ for all } x\in U.
\tag{$\dagger$}
\]

For each $i$, choose $\chi_i\in C_c(\C)$ with $0\le \chi_i\le 1$, $\supp(\chi_i)\subset D_i$,
and $\chi_i(a_i)=1$. Then $\langle F(x_0),\chi_i\rangle=1$, so by continuity there is a neighborhood
$U_i\subset U$ of $x_0$ such that $\langle F(x),\chi_i\rangle>1/2$ for all $x\in U_i$.
In particular, $F(x)\cap D_i\neq\varnothing$ for $x\in U_i$ (otherwise the sum is $0$).
Let
\[
U':=\bigcap_{i=1}^n U_i\subset U.
\]
Then for all $x\in U'$:
\begin{itemize}
\item $F(x)\cap D_i\neq\varnothing$ for each $i$;
\item by ($\dagger$), every point of $F(x)\cap B^-$ lies in $\bigcup_i D_i$.
\end{itemize}
Since the disks are disjoint and $\#(F(x)\cap B^-)=n$, it follows that for each $x\in U'$
\[
\#(F(x)\cap D_i)=1\quad\text{for all }i,
\qquad\text{and}\qquad
F(x)\cap B^-=\bigsqcup_{i=1}^n \bigl(F(x)\cap D_i\bigr).
\tag{$\ddagger$}
\]
Hence we may define, for $x\in U'$, the unique point
\[
a_i(x)\in D_i\quad\text{by}\quad F(x)\cap D_i=\{a_i(x)\}.
\]
Define
\[
\widetilde r(x):=\bigl(a_1(x),\dots,a_n(x)\bigr)\in D_1\times\cdots\times D_n.
\]

\medskip
\noindent\textbf{Claim 1.} Each map $a_i:U'\to D_i$ is continuous.

\emph{Proof.}
Let $x_\alpha\to x$ be a net in $U'$. Put $p_\alpha:=a_i(x_\alpha)\in \overline{D_i}$.
Since $\overline{D_i}$ is compact, $(p_\alpha)$ has a convergent subnet $p_{\alpha_\beta}\to p\in \overline{D_i}$.
We show $p=a_i(x)$, which implies the whole net converges and hence continuity.

Choose an open neighborhood $W$ of $p$ with $\overline{W}\subset D_i$. Let $\psi\in C_c(\C)$ satisfy
$0\le\psi\le1$, $\psi\equiv 1$ on $\overline{W}$, and $\supp(\psi)\subset D_i$.
Then for $\beta$ large, $p_{\alpha_\beta}\in W$, hence
\[
\langle F(x_{\alpha_\beta}),\psi\rangle
=\sum_{q\in F(x_{\alpha_\beta})}\psi(q)
=\psi(p_{\alpha_\beta})=1
\]
using ($\ddagger$) (there is exactly one point in $D_i$).
By vague continuity, $\langle F(x_{\alpha_\beta}),\psi\rangle\to \langle F(x),\psi\rangle$, so
$\langle F(x),\psi\rangle=1$. Since $\supp(\psi)\subset D_i$ and $F(x)\cap D_i=\{a_i(x)\}$ by ($\ddagger$),
we have $\langle F(x),\psi\rangle=\psi(a_i(x))$, hence $\psi(a_i(x))=1$ and therefore $a_i(x)\in \overline{W}$.
As $W$ was an arbitrary neighborhood of $p$ with $\overline{W}\subset D_i$, this forces $p=a_i(x)$.
\qed

\medskip
By Claim~1, $\widetilde r:U'\to D_1\times\cdots\times D_n$ is continuous.
Let
\[
q:D_1\times\cdots\times D_n \longrightarrow C_n(\operatorname{int}B),
\qquad
q(z_1,\dots,z_n):=\{z_1,\dots,z_n\}.
\]
Because the disks are pairwise disjoint, $q$ is a homeomorphism onto the open subset of
$C_n(\operatorname{int}B)$ consisting of configurations having exactly one point in each $D_i$
(a standard configuration-space chart).
For $x\in U'$, we have by ($\ddagger$)
\[
r_B(x)=F(x)\cap \operatorname{int}B = F(x)\cap B^- = \{a_1(x),\dots,a_n(x)\}= q(\widetilde r(x)).
\]
Hence $r_B|_{U'}=q\circ \widetilde r$ is continuous at $x_0$.
Since $x_0$ was arbitrary in $K_0$, $r_B$ is continuous on $K_0$.
\end{proof}

%------------------------------------------------------------
% Auxiliary lemmas: uniform room + uniform separation + uniform time partition
%------------------------------------------------------------

\begin{lemma}[Uniform room inside $U$ and uniform separation]\label{lem:room-sep}
Let $K$ be a finite CW complex, $P$ a finite discrete space, and
$\eta:K\times[0,1]\times P\to \C$ a continuous map such that for each $(x,t)$
the map $\eta_{x,t}:P\to\C$ is injective.
Let $U\subset\C$ be open and assume $\eta(K\times[0,1]\times P)\subset U$.
Then there exist constants $\delta>0$ and $\sigma>0$ such that:
\begin{enumerate}
\item (\emph{Room}) If $c\in \eta(K\times[0,1]\times P)$ then $\overline{B(c,\delta)}\subset U$.
\item (\emph{Separation}) For all $(x,t)\in K\times[0,1]$ and all distinct $p\neq q\in P$,
      \[
      |\eta(x,t,p)-\eta(x,t,q)|\ge \sigma.
      \]
\end{enumerate}
\end{lemma}

\begin{proof}
Set $B:=K\times[0,1]\times P$. Since $K$ is compact (finite CW) and $P$ is finite,
$B$ is compact. Hence $C:=\eta(B)\subset U$ is compact.

(1) Since $U$ is open and $C\Subset U$, the distance
\[
\delta:=\dist(C,\C\setminus U)=\inf\{|c-y|:c\in C,\,y\in\C\setminus U\}
\]
is strictly positive. Then $\overline{B(c,\delta)}\subset U$ for all $c\in C$.

(2) For each ordered pair $(p,q)$ with $p\neq q$, define
\[
d_{p,q}:K\times[0,1]\to(0,\infty),\qquad d_{p,q}(x,t):=|\eta(x,t,p)-\eta(x,t,q)|.
\]
Each $d_{p,q}$ is continuous and positive everywhere (injectivity).
By compactness of $K\times[0,1]$, each $d_{p,q}$ attains a positive minimum.
Since there are finitely many pairs $(p,q)$, the finite minimum
\[
\sigma:=\min_{p\neq q}\ \min_{(x,t)\in K\times[0,1]} d_{p,q}(x,t)
\]
is positive, and the stated inequality follows.
\end{proof}

\begin{lemma}[Uniform time subdivision]\label{lem:uniform-partition}
In the setup of Lemma~\ref{lem:room-sep}, fix $r>0$.
Then there exists a partition $0=t_0<t_1<\cdots<t_N=1$ such that for every $i$,
every $x\in K$, every $p\in P$, and every $t\in[t_{i-1},t_i]$,
\[
|\eta(x,t,p)-\eta(x,t_{i-1},p)|<\frac r6.
\]
\end{lemma}

\begin{proof}
The map $\eta$ is uniformly continuous on the compact domain $K\times[0,1]\times P$,
so there exists $\Delta>0$ such that $|t-s|<\Delta$ implies
$|\eta(x,t,p)-\eta(x,s,p)|<r/6$ for all $(x,p)\in K\times P$.
Choose a partition with mesh $<\Delta$.
\end{proof}

%------------------------------------------------------------
% Simultaneous finite-point bump move on one slab
%------------------------------------------------------------

\begin{lemma}[Simultaneous bump move on a time slab]\label{lem:simul-bump}
Let $m\in\N$, $K$ be compact, and $U\subset\C$ be open.
Fix numbers $r>0$ and $\sigma>0$ such that $4r/3<\sigma$.
Assume we are given continuous maps
\[
a_j:K\to U,\qquad b_j:K\times[t_{-},t_{+}]\to U\qquad (j=1,\dots,m)
\]
such that
\begin{enumerate}
\item (\emph{small displacement}) $|b_j(x,t)-a_j(x)|<r/6$ for all $(x,t)$ and all $j$;
\item (\emph{separated centers}) $|a_j(x)-a_k(x)|\ge \sigma$ for all $x\in K$ and all $j\neq k$;
\item (\emph{room}) $\overline{B(a_j(x),2r/3)}\subset U$ for all $x\in K$ and all $j$.
\end{enumerate}
Then there exists a continuous map
\[
H:K\times[t_{-},t_{+}]\times \C\to\C,\qquad H_{x,t}:=H(x,t,\cdot),
\]
such that for each $(x,t)$:
\begin{enumerate}
\item $H_{x,t}:\C\to\C$ is a homeomorphism;
\item $H_{x,t}=\id$ on $\C\setminus U$;
\item $H_{x,t}(a_j(x))=b_j(x,t)$ for every $j$;
\item if additionally $b_j(x,t_{-})=a_j(x)$ for all $j$, then $H_{x,t_{-}}=\id$.
\end{enumerate}
\end{lemma}

\begin{proof}
Fix once and for all a Lipschitz bump function $\theta:[0,\infty)\to[0,1]$
as in Lemma~\ref{lem:bumpmove} for the radius $r$:
$\theta\equiv 1$ on $[0,r/3]$ and $\theta\equiv 0$ on $[2r/3,\infty)$.

For each $j$ define
\[
H_j(x,t)(z):=z+\bigl(b_j(x,t)-a_j(x)\bigr)\,\theta\bigl(|z-a_j(x)|\bigr).
\]
By Lemma~\ref{lem:bumpmove} (applied with center $a=a_j(x)$ and target point $p=b_j(x,t)$),
each $H_j(x,t)$ is a homeomorphism, satisfies $H_j(x,t)(a_j(x))=b_j(x,t)$,
and is the identity outside $B(a_j(x),2r/3)$.
The map $(x,t,z)\mapsto H_j(x,t)(z)$ is continuous because $a_j,b_j$ are continuous and
$(p,z)\mapsto h_{a,p}(z)$ is continuous in Lemma~\ref{lem:bumpmove}.

By the separation hypothesis and $4r/3<\sigma$, the balls
$B(a_j(x),2r/3)$ are pairwise disjoint for fixed $x$. Hence the supports of
$H_j(x,t)$ are disjoint, so the $H_j(x,t)$ commute.

Define the simultaneous move
\[
H_{x,t}:=H_m(x,t)\circ\cdots\circ H_1(x,t).
\]
This is a homeomorphism (finite composition of homeomorphisms), and it equals $\id$
on $\C\setminus U$ since each $H_j$ does and each support lies in $U$ by (3).
Moreover, for each $j$, the point $a_j(x)$ lies in the support of $H_j$ but in no
support of $H_k$ for $k\neq j$ (disjointness), hence
\[
H_{x,t}(a_j(x))=H_j(x,t)(a_j(x))=b_j(x,t).
\]
If $b_j(x,t_-)=a_j(x)$ for all $j$, then each $H_j(x,t_-)=\id$ and hence $H_{x,t_-}=\id$.
Continuity of $(x,t,z)\mapsto H_{x,t}(z)$ follows from continuity of each $H_j$ and
finite composition.
\end{proof}

%------------------------------------------------------------
% Main theorem for M = C
%------------------------------------------------------------

The following result is inspired by the isotopy extension theorem of Edwards-–Kirby, 
together with their relative and support-control deformation theory (Theorem~5.1 and the 
relative discussion in section 7; see~\cite{EK71}).
We also benefit from the discussion surrounding \cite[Theorem~3.9]{Kup15}.

\begin{theorem}[Parametric Edwards-Kirby isotopy extension]\label{thm:EK}
Let $U\subset\C$ be open, and let $K$ be a finite CW complex with subcomplex $L\subset K$.
Let $P$ be a finite discrete space.
Suppose $\eta:K\times[0,1]\times P\to U$ is continuous and:
\begin{enumerate}
\item for each $(x,t)$, $\eta_{x,t}:P\to U$ is injective;
\item for every $x\in L$, $\eta_{x,t}=\eta_{x,0}$ for all $t\in[0,1]$.
\end{enumerate}
Then there exists a continuous map
\[
\Phi:K\times[0,1]\times\C\to\C,\qquad \Phi_{x,t}:=\Phi(x,t,\cdot),
\]
such that:
\begin{enumerate}
\item for each $(x,t)$, $\Phi_{x,t}$ is a homeomorphism of $\C$;
\item $\Phi_{x,0}=\id$ for all $x$;
\item $\Phi_{x,t}=\id$ on $\C\setminus U$ for all $(x,t)$;
\item for all $(x,t,p)$, $\Phi_{x,t}(\eta(x,0,p))=\eta(x,t,p)$;
\item for all $x\in L$ and all $t$, $\Phi_{x,t}=\id$.
\end{enumerate}
Moreover, if for some subinterval $I\subset[0,1]$ the map $t\mapsto \eta(x,t,p)$ is constant on $I$
for all $p\in P$ (uniformly in $x$), then $\Phi_{x,t}$ can be chosen constant on $I$ as well.
\end{theorem}

\begin{proof}
Let $P=\{p_1,\dots,p_m\}$.

\smallskip\noindent
\textbf{(A) Choose uniform radii.}
Apply Lemma~\ref{lem:room-sep} to obtain $\delta>0$ and $\sigma>0$ such that:
(i) $\overline{B(c,\delta)}\subset U$ for all $c$ in the compact set $C:=\eta(K\times[0,1]\times P)$,
and (ii) $|\eta(x,t,p_j)-\eta(x,t,p_k)|\ge\sigma$ for all $j\neq k$ and all $(x,t)$.
Choose $r>0$ such that
\[
\frac{2r}{3}<\delta,\qquad \frac{4r}{3}<\sigma.
\tag{$*$}
\]

\smallskip\noindent
\textbf{(B) Choose a uniform time subdivision.}
Apply Lemma~\ref{lem:uniform-partition} to obtain a partition
$0=t_0<t_1<\cdots<t_N=1$ such that for every $i$, every $x\in K$, every $j$,
and every $t\in[t_{i-1},t_i]$,
\[
|\eta(x,t,p_j)-\eta(x,t_{i-1},p_j)|<\frac r6.
\tag{$**$}
\]

\smallskip\noindent
\textbf{(C) Build the ambient isotopy slab by slab.}
For each $i=1,\dots,N$ and each $j=1,\dots,m$, define
\[
a_{i,j}(x):=\eta(x,t_{i-1},p_j),\qquad b_{i,j}(x,t):=\eta(x,t,p_j)\quad (t\in[t_{i-1},t_i]).
\]
Then $(**)$ gives the small displacement condition for Lemma~\ref{lem:simul-bump}.
The separation condition for the centers $a_{i,j}(x)$ follows from the uniform separation $\sigma$,
and the room condition $\overline{B(a_{i,j}(x),2r/3)}\subset U$ follows from $(*)$ because each
$a_{i,j}(x)\in C$ and $2r/3<\delta$.

Hence Lemma~\ref{lem:simul-bump} produces, for each $i$, a continuous family of
homeomorphisms
\[
H_i:K\times[t_{i-1},t_i]\times\C\to\C,\qquad (x,t,z)\mapsto H_i(x,t)(z),
\]
such that for each $(x,t)\in K\times[t_{i-1},t_i]$:
\begin{enumerate}
\item $H_i(x,t)=\id$ on $\C\setminus U$;
\item $H_i(x,t)\bigl(\eta(x,t_{i-1},p_j)\bigr)=\eta(x,t,p_j)$ for all $j$;
\item $H_i(x,t_{i-1})=\id$.
\end{enumerate}

Now define $\Phi_{x,t}$ inductively in $i$.
Set $\Phi_{x,t_0}:=\id$. Suppose $\Phi_{x,t_{i-1}}$ is defined.
For $t\in[t_{i-1},t_i]$, define
\[
\Phi_{x,t}:=H_i(x,t)\circ \Phi_{x,t_{i-1}}.
\tag{$\dagger$}
\]
Then each $\Phi_{x,t}$ is a homeomorphism, and $\Phi_{x,t}=\id$ on $\C\setminus U$
since both factors in $(\dagger)$ are.

\smallskip\noindent
\textbf{(D) Continuity.}
We claim $(x,t,z)\mapsto \Phi_{x,t}(z)$ is continuous on each slab and glues across endpoints.
Assume by induction that $(x,z)\mapsto \Phi_{x,t_{i-1}}(z)$ is continuous.
On $t\in[t_{i-1},t_i]$ we have
\[
\Phi(x,t,z)=H_i\bigl(x,t,\Phi(x,t_{i-1},z)\bigr),
\]
a composition of continuous maps, hence continuous.
At $t=t_{i-1}$, since $H_i(x,t_{i-1})=\id$, the formula agrees with the already defined
$\Phi_{x,t_{i-1}}$, so the pieces glue continuously. Hence $\Phi$ is continuous on
$K\times[0,1]\times\C$.

\smallskip\noindent
\textbf{(E) Extension property.}
Fix $j$ and proceed by induction on $i$.
For $i=0$, $\Phi_{x,t_0}(\eta(x,0,p_j))=\eta(x,0,p_j)$.
Assume $\Phi_{x,t_{i-1}}(\eta(x,0,p_j))=\eta(x,t_{i-1},p_j)$.
Then for $t\in[t_{i-1},t_i]$,
\[
\Phi_{x,t}(\eta(x,0,p_j))
=H_i(x,t)\bigl(\Phi_{x,t_{i-1}}(\eta(x,0,p_j))\bigr)
=H_i(x,t)\bigl(\eta(x,t_{i-1},p_j)\bigr)
=\eta(x,t,p_j).
\]
This proves (4).

\smallskip\noindent
\textbf{(F) Initial time and support control.}
By construction $\Phi_{x,0}=\id$.
Support control holds since each $H_i(x,t)=\id$ on $\C\setminus U$ and the recursion
preserves this property.

\smallskip\noindent
\textbf{(G) Relative property on $L$.}
If $x\in L$, then by hypothesis $\eta(x,t,p_j)=\eta(x,0,p_j)$ for all $t$ and all $j$.
In particular, on each slab $b_{i,j}(x,t)\equiv a_{i,j}(x)$, so Lemma~\ref{lem:simul-bump}
gives $H_i(x,t)\equiv\id$, and then $(\dagger)$ yields $\Phi_{x,t}\equiv\id$ for all $t$.

\smallskip\noindent
\textbf{(H) Optional stationarity in $t$.}
If $\eta(x,t,p)$ is constant in $t$ on a subinterval $I$, refine the partition so that
$I$ is a union of whole slabs. On those slabs, $b_{i,j}(x,t)\equiv a_{i,j}(x)$, hence
$H_i(x,t)\equiv\id$ and $(\dagger)$ shows $\Phi_{x,t}$ is constant on $I$.
\end{proof}

\subsection{Infinite gap lemma}
For the following result, see \cite{Engelking}.
\begin{theorem}[Katětov--Tong insertion theorem]
\label{thm:KT}
Let $X$ be a normal topological space. If
$a:X\to\R$ is upper semicontinuous, $b:X\to\R$ is lower semicontinuous,
and $a(x)<b(x)$ for all $x$, then there exists a continuous $c:X\to\R$
with $a(x)<c(x)<b(x)$ for all $x$.
\end{theorem}

%%%%%%%%%%%%%%%%%%%%%%%%%%%%%%%%%%%%%%%%%%%%%%%%%%%%%%%%%%%%%%%%%%%%%%%%%%%%%%
%  Correct boundary normalization lemma (patched so chosen radii are OFF boundary)
%%%%%%%%%%%%%%%%%%%%%%%%%%%%%%%%%%%%%%%%%%%%%%%%%%%%%%%%%%%%%%%%%%%%%%%%%%%%%%

\begin{lemma}[Boundary normalization near two circles]
\label{lem:boundary-normalize-correct}
Fix radii $0<R<2R<3R$ and put
\[
A:=\{R\le |z|\le 2R\},\qquad B:=\{R/2\le |z|\le 3R\}.
\]
Fix $\varepsilon>0$ so small that
\[
R/2+\varepsilon<R
\qquad\text{and}\qquad
2R<3R-\varepsilon.
\]
Let $K$ be compact metric space and $F:K\to C^{lf}_\infty(\C)$ be continuous where $C^{lf}_{\infty}(\C)$ is endowed with the vague topology.

Define the two one-sided collar functions
\[
u_{\mathrm{in}}(x):=\sup\{|z|:\ z\in F(x),\ R/2-\varepsilon\le |z|<R/2\},
\quad (\sup\emptyset:=R/2-\varepsilon),
\]
\[
v_{\mathrm{out}}(x):=\inf\{|z|:\ z\in F(x),\ 3R< |z|\le 3R+\varepsilon\},
\quad (\inf\emptyset:=3R+\varepsilon).
\]
Then $u_{\mathrm{in}}$ is upper semicontinuous and satisfies $u_{\mathrm{in}}(x)<R/2$ for all $x$,
and $v_{\mathrm{out}}$ is lower semicontinuous and satisfies $3R<v_{\mathrm{out}}(x)$ for all $x$.

Consequently, by Theorem~\ref{thm:KT} there exist continuous functions
$r^{\mathrm{in}},r^{\mathrm{out}}:K\to\R$ such that for all $x$,
\[
u_{\mathrm{in}}(x)<r^{\mathrm{in}}(x)<R/2,
\qquad
3R<r^{\mathrm{out}}(x)<v_{\mathrm{out}}(x).
\]
In particular,
\[
F(x)\cap\{|z|=r^{\mathrm{in}}(x)\}=\varnothing,\qquad
F(x)\cap\{|z|=r^{\mathrm{out}}(x)\}=\varnothing
\quad \mbox{ for all } x.
\]

Moreover, there is a continuous family of radial homeomorphisms
$\Psi_{x,t}:\C\to\C$ ($x\in K$, $t\in[0,1]$) such that:
\begin{enumerate}
\item $\Psi_{x,0}=\mathrm{id}$;
\item $\Psi_{x,t}=\mathrm{id}$ on $\{|z|\le R/2-\varepsilon\}\cup A\cup\{|z|\ge 3R+\varepsilon\}$;
\item each $\Psi_{x,t}$ is supported in the fixed compact set
\[
S:=\{R/2-\varepsilon\le |z|\le R/2+\varepsilon\}\ \cup\ \{3R-\varepsilon\le |z|\le 3R+\varepsilon\};
\]
\item $\Psi_{x,1}$ satisfies
\[
\Psi_{x,1}(\{|z|=r^{\mathrm{in}}(x)\})=\{|z|=R/2\},
\qquad
\Psi_{x,1}(\{|z|=r^{\mathrm{out}}(x)\})=\{|z|=3R\};
\]
\item consequently,
\[
(\Psi_{x,1})_*(F(x))\cap \partial B=\varnothing\qquad\mbox{ for all } x.
\]
\end{enumerate}
\end{lemma}

\begin{proof}

\smallskip
\textbf{(A) Upper semicontinuity of $u_{\mathrm{in}}$ and $u_{\mathrm{in}}(x)<R/2$.}
Fix $c<R/2$. For $n\ge1$ set
\[
Q_{c,n}^{\mathrm{in}}:=\{z:\ c\le |z|\le R/2-1/n\}\cap\{R/2-\varepsilon\le |z|\le R/2\},
\]
compact. Then
\[
\{x:\ u_{\mathrm{in}}(x)<c\}
=
\bigcap_{n\ge1}\{x:\ F(x)\cap Q_{c,n}^{\mathrm{in}}=\varnothing\},
\]
so it is open by Lemma \ref{lem:openness}. Hence $u_{\mathrm{in}}$ is upper semicontinuous.
Moreover, $F(x)\cap\{R/2-\varepsilon\le |z|\le R/2\}$ is finite and all its radii are $<R/2$,
so the supremum is $<R/2$ (or equals $R/2-\varepsilon$ if empty).

\smallskip
\textbf{(B) Lower semicontinuity of $v_{\mathrm{out}}$ and $3R<v_{\mathrm{out}}(x)$.}
Fix $c>3R$. For $n\ge1$ set
\[
P_{c,n}^{\mathrm{out}}:=\{z:\ 3R+1/n\le |z|\le c\}\cap\{3R\le |z|\le 3R+\varepsilon\},
\]
compact. Then
\[
\{x:\ v_{\mathrm{out}}(x)>c\}
=
\bigcap_{n\ge1}\{x:\ F(x)\cap P_{c,n}^{\mathrm{out}}=\varnothing\},
\]
is open by Lemma \ref{lem:openness}. Hence $v_{\mathrm{out}}$ is lower semicontinuous.
Also all radii in the defining set for $v_{\mathrm{out}}$ are $>3R$, so $v_{\mathrm{out}}(x)>3R$.

\smallskip
\textbf{(C) Choose continuous empty radii off the boundary.}
Apply Theorem~\ref{thm:KT} to $a(x)=u_{\mathrm{in}}(x)$ and $b(x)=R/2$ to obtain
continuous $r^{\mathrm{in}}$ with $u_{\mathrm{in}}<r^{\mathrm{in}}<R/2$.
Apply Theorem~\ref{thm:KT} to $a(x)=3R$ and $b(x)=v_{\mathrm{out}}(x)$ to obtain
continuous $r^{\mathrm{out}}$ with $3R<r^{\mathrm{out}}<v_{\mathrm{out}}$.

Then $\{|z|=r^{\mathrm{in}}(x)\}$ is empty in $F(x)$: any point of that radius would lie in the collar below $R/2$
and contradict $r^{\mathrm{in}}(x)>u_{\mathrm{in}}(x)$. Similarly $\{|z|=r^{\mathrm{out}}(x)\}$ is empty:
any point of that radius would lie in the collar above $3R$ and contradict $r^{\mathrm{out}}(x)<v_{\mathrm{out}}(x)$.

\smallskip
\textbf{(D) Radial homeomorphisms with fixed support and identity on $A$.}
For each $x$, define a strictly increasing homeomorphism $\sigma_x:[0,\infty)\to[0,\infty)$ by:
\[
\sigma_x(r)=r\ \text{for }r\le R/2-\varepsilon,\ \text{for }R\le r\le 2R,\ \text{and for }r\ge 3R+\varepsilon,
\]
and let $\sigma_x$ be piecewise linear on $[R/2-\varepsilon,R/2+\varepsilon]$ and on $[3R-\varepsilon,3R+\varepsilon]$
with
\[
\sigma_x(R/2-\varepsilon)=R/2-\varepsilon,\
\sigma_x(r^{\mathrm{in}}(x))=R/2,\
\sigma_x(R/2+\varepsilon)=R/2+\varepsilon,
\]
\[
\sigma_x(3R-\varepsilon)=3R-\varepsilon,\
\sigma_x(r^{\mathrm{out}}(x))=3R,\
\sigma_x(3R+\varepsilon)=3R+\varepsilon.
\]
Define $\Psi_x(0)=0$ and $\Psi_x(re^{i\theta})=\sigma_x(r)e^{i\theta}$ for $r>0$.
Then $\Psi_x$ is a homeomorphism of $\C$, supported in the fixed compact set $S$, and equals the identity on $A$.

Define the isotopy
\[
\sigma_{x,t}(r)=(1-t)r+t\,\sigma_x(r),\qquad
\Psi_{x,t}(0)=0,\ \Psi_{x,t}(re^{i\theta})=\sigma_{x,t}(r)e^{i\theta}.
\]
Each $\sigma_{x,t}$ is strictly increasing, hence each $\Psi_{x,t}$ is a homeomorphism, and
$(x,t,z)\mapsto\Psi_{x,t}(z)$ is continuous.

\smallskip
\textbf{(E) Boundary avoidance.}
Since $\sigma_x$ is strictly increasing,
\[
\Psi_x^{-1}(\{|z|=R/2\})=\{|z|=r^{\mathrm{in}}(x)\},
\qquad
\Psi_x^{-1}(\{|z|=3R\})=\{|z|=r^{\mathrm{out}}(x)\}.
\]
Both preimage circles are empty in $F(x)$, hence $(\Psi_x)_*(F(x))\cap\partial B=\varnothing$.
\end{proof}

The following is a standard result in algebraic topology.
\begin{lemma}\label{lem:sc-to-kpi1}
Let $K$ be a simply connected CW complex and let $Y$ be a $K(\pi,1)$.
Then every continuous map $f:K\to Y$ is null-homotopic.
In particular, if $U\subset\C$ is an open disk or open annulus, for $n\in \N$,
every continuous map $f:K\to Conf_n(U)$
is null-homotopic.
\end{lemma}

%%%%%%%%%%%%%%%%%%%%%%%%%%%%%%%%%%%%%%%%%%%%%%%%%%%%%%%%%%%%%%%%%%%%%%%%%%%%%%
%  Infinite Gap Lemma (final)
%%%%%%%%%%%%%%%%%%%%%%%%%%%%%%%%%%%%%%%%%%%%%%%%%%%%%%%%%%%%%%%%%%%%%%%%%%%%%%
The following result is the key technical tool we use to prove our main results.

\begin{lemma}[Infinite Gap Lemma]
\label{lem:inf-gap}
Let $K$ be a finite CW complex with $\pi_1(K)=0$ and let
$F:K\to C^{lf}_\infty(\C)$ be continuous (vague topology).
Then there exist radii $0<R_1<R_2<\cdots\to\infty$ and a continuous map
$G:K\to C^{lf}_\infty(\C)$ homotopic to $F$ such that for every $m\ge1$ and every $x\in K$,
\[
G(x)\cap A_m=\varnothing \qquad \mbox{ where } \qquad A_m:=\{z\in\C:\ R_m\le |z|\le 2R_m\}.
\]
Moreover, the homotopy $F\simeq G$ is a countable concatenation of homotopies,
where stage $m$ is induced by composing with a $K$--parametrized compactly supported
ambient isotopy supported in a compact annulus, and these supports may be chosen pairwise disjoint.
\end{lemma}

\begin{proof}
\noindent\textbf{Step 0 (choose disjoint supports).}
Choose radii $R_m$ inductively so that $R_{m+1}>10^4R_m$.
Put
\[
B_m:=\{R_m/2\le |z|\le 3R_m\},\quad U_m:=\operatorname{int}(B_m),\quad
A_m:=\{R_m\le |z|\le 2R_m\}\subset U_m.
\]
Fix $\delta_m:=R_m/100$ and
\[
W_m:=\{R_m/2-\delta_m<|z|<3R_m+\delta_m\}.
\]
Then $\overline{W_m}$ are pairwise disjoint compact sets.

We inductively construct continuous maps $F^{(m)}:K\to C^{lf}_\infty(\C)$ with $F^{(0)}:=F$ such that
\[
F^{(m)}(x)\cap A_j=\varnothing\quad\text{for all }j\le m\text{ and all }x,
\]
and $F^{(m-1)}\simeq F^{(m)}$ by a homotopy supported in $\overline{W_m}$.

\medskip
\noindent\textbf{Stage $m$ (clear $A_m$).}
Assume $F^{(m-1)}$ is constructed.

\smallskip
\emph{(1) Normalize boundaries.}
Choose $\varepsilon_m>0$ with $\varepsilon_m<\delta_m$ and satisfying
$R_m/2+\varepsilon_m<R_m$ and $2R_m<3R_m-\varepsilon_m$.
Apply Lemma~\ref{lem:boundary-normalize-correct} (with $R=R_m$, $\varepsilon=\varepsilon_m$) to
$F^{(m-1)}:K\to C^{lf}_\infty(\C)$.
We obtain a continuous isotopy $\Psi^{(m)}_{x,t}$ supported in a fixed compact subset
$S_m\subset W_m$ and equal to the identity on $A_m$, such that
\[
\widetilde F^{(m-1)}(x):=(\Psi^{(m)}_{x,1})_*\bigl(F^{(m-1)}(x)\bigr)
\quad\text{satisfies}\quad
\widetilde F^{(m-1)}(x)\cap\partial B_m=\varnothing\ \ \mbox{ for all } x.
\]
By Lemma~\ref{lem:param-push}, $x\mapsto \widetilde F^{(m-1)}(x)$ is continuous, and this step yields a homotopy
$F^{(m-1)}\simeq \widetilde F^{(m-1)}$ supported in $\overline{W_m}$. Since $\Psi^{(m)}_{x,t}$ is the identity on $A_m$,
this step does not change membership in $A_m$.

\smallskip
\emph{(2) Restrict to $U_m$ with constant cardinality.}
Apply Lemma~\ref{lem:restrict} to $B=B_m$ and the map $\widetilde F^{(m-1)}$.
On each connected component of $K$, the integer
\[
n_m:=\#\bigl(\widetilde F^{(m-1)}(x)\cap B_m\bigr)
\]
is constant, and the restriction map
\[
r_m:K\to C_{n_m}(U_m),\qquad r_m(x):=\widetilde F^{(m-1)}(x)\cap U_m
\]
is continuous.

\smallskip
\emph{(3) Lift and null-homotope into $U_m\setminus A_m$.}
Let
\[
V_m:=\{(5/2)R_m<|z|<3R_m\}\subset U_m\setminus A_m.
\]
Since $\pi_1(K)=0$ and $Conf_{n_m}(U_m)\to C_{n_m}(U_m)$ is a covering, $r_m$ admits a continuous lift
$\bar r_m:K\to Conf_{n_m}(U_m)$.
By the fact that $Conf_{n_m}(U_m)$ is a $K(\pi,1)$ and Lemma~\ref{lem:sc-to-kpi1},
$\bar r_m$ is null-homotopic. Fix $b_m\in Conf_{n_m}(V_m)$ and choose a homotopy
\[
h_m:K\times[0,1]\to Conf_{n_m}(U_m)
\]
with $h_m(\cdot,0)=\bar r_m$ and $h_m(\cdot,1)\equiv b_m$.

\smallskip
\emph{(4) Realize $h_m$ by EK extension, supported in $U_m$.}
Let $P_m=\{1,\dots,n_m\}$ and define
\[
\eta_m:K\times[0,1]\times P_m\to U_m,\qquad \eta_m(x,t,i)=(h_m(x,t))_i.
\]
By Theorem~\ref{thm:EK} (Edwards--Kirby) with support control, there exists a continuous family of
homeomorphisms $\Phi^{(m)}_{x,t}:\C\to\C$ such that:
\begin{itemize}
\item $\Phi^{(m)}_{x,0}=\mathrm{id}$;
\item $\Phi^{(m)}_{x,t}=\mathrm{id}$ on $\C\setminus U_m$ (hence supported in $B_m\subset W_m$);
\item $(\Phi^{(m)}_{x,t})_*$ carries $\bar r_m(x)$ to $h_m(x,t)$.
\end{itemize}
Define
\[
F^{(m)}(x):=(\Phi^{(m)}_{x,1})_*\bigl(\widetilde F^{(m-1)}(x)\bigr).
\]
Then $F^{(m)}(x)\cap U_m$ is the unordered configuration underlying $b_m\subset V_m\subset U_m\setminus A_m$,
so
\[
F^{(m)}(x)\cap A_m=\varnothing\qquad\mbox{ for all } x.
\]
By Lemma~\ref{lem:param-push}, the stage homotopy is continuous in the vague topology and supported in $\overline{W_m}$.

\smallskip
\emph{(5) Earlier gaps remain empty.}
For $j<m$, we have $A_j\subset\{|z|\le 2R_{m-1}\}$, while
$\overline{W_m}\subset\{|z|\ge R_m/2-\delta_m\}$ and $R_m/2-\delta_m\gg 2R_{m-1}$ by the choice $R_m>10^4R_{m-1}$.
Hence $A_j\cap \overline{W_m}=\varnothing$. Therefore stage $m$ does not affect earlier gaps. Since the $\overline{W_m}$
are pairwise disjoint, later stages do not affect earlier ones either.

\medskip
\noindent\textbf{Step 1 (define $G$ by stabilization on compacta).}
Define $G(x)\subset\C$ by prescribing its intersection with every compact set $C$:
choose any $M=M(C)$ with $C\cap\overline{W_m}=\varnothing$ for all $m\ge M$, and set
\[
G(x)\cap C \;:=\; F^{(M)}(x)\cap C.
\]
This is well-defined: if $M'\ge M$ is another such index, then by stabilization,
\[
F^{(M')}(x)\cap C = F^{(M)}(x)\cap C,
\]
so the right-hand side does not depend on the choice of $M$.

Equivalently, since the intersections stabilize, for any sufficiently large $M$ one may write
\[
G(x)\cap C \;=\; \bigcap_{m\ge M}\bigl(F^{(m)}(x)\cap C\bigr)
\;=\; F^{(M)}(x)\cap C .
\]

Then $G(x)\in C^{lf}_\infty(\C)$ and $G(x)\cap A_m=\varnothing$ for all $m$.

\smallskip
\noindent\emph{Continuity of $G$.}
For $\varphi\in C_c(\C)$ let $C=\supp(\varphi)$ and choose $M$ as above. Then
\[
\langle G(x),\varphi\rangle=\langle F^{(M)}(x),\varphi\rangle,
\]
which is continuous in $x$, hence $G$ is continuous in the vague topology.

\smallskip
\noindent\emph{Homotopy $F\simeq G$.}
Concatenate the stage homotopies using $t_m=1-2^{-m}$. Vague continuity is checked by test functions:
for fixed $\varphi\in C_c(\C)$, the compact $\supp(\varphi)$ meets only finitely many $\overline{W_m}$, so only finitely many
stages affect $\langle\cdot,\varphi\rangle$, and on each time slab the evaluation is continuous. Thus the concatenation
defines a continuous homotopy $F\simeq G$.
\end{proof}

The following result is straight forward.
\begin{lemma}[Basepoint change for $m\ge2$]\label{lem:basepoint-change}
Let $X$ be path-connected, $m\ge2$, and $\gamma:[0,1]\to X$ a path from $x_0$ to $x_1$.
Then there is a natural isomorphism $\gamma_{\#}:\pi_m(X,x_0)\cong \pi_m(X,x_1)$.
\end{lemma}

\begin{theorem}[Asphericity of $C^{lf}_\infty(\C)$]
\label{thm:aspherical-corrected}
The space $C^{lf}_{\infty}(\C)$ is aspherical.
\end{theorem}

\begin{proof}
Fix $k\ge2$ and a based map
\[
f:(S^k,*)\longrightarrow \bigl(C^{lf}_\infty(\C), \N\bigr).
\]

\medskip
\noindent\textbf{Step 0 (gap lemma with geometric growth).}
Apply the Infinite Gap Lemma to obtain a continuous map
\[
g:S^k\to C^{lf}_\infty(\C)
\]
homotopic to $f$ and radii $0<R_1<R_2<\cdots\to\infty$ such that
\[
g(x)\cap A_m=\varnothing,
\qquad
A_m:=\{z\in\C:\ R_m\le |z|\le 2R_m\},
\quad \mbox{ for all } x\in S^k,\ \mbox{ for all } m\ge1.
\]
By passing to a subsequence (and relabeling), we may assume in addition that
\begin{equation}
\label{eq:Rm-growth}
R_m>2R_{m-1}\qquad(m\ge2).
\end{equation}
(Indeed, since $R_m\to\infty$, choose inductively $m_1<m_2<\cdots$ with
$R_{m_{\ell+1}}>2R_{m_\ell}$; then $g$ avoids the annuli $A_{m_\ell}$ as well.)
Fix $Q:=g(*)$.

\medskip
\noindent\textbf{Step 1 (shell decomposition separated by empty annuli).}
Define open shells
\[
S_1:=\{|z|<R_1\},
\qquad
S_m:=\{2R_{m-1}<|z|<R_m\}\quad(m\ge2).
\]
Condition \eqref{eq:Rm-growth} guarantees $S_m\neq\varnothing$ for all $m\ge2$.
Moreover,
\[
\partial\overline{S_1}\subset \{|z|=R_1\}\subset A_1,
\qquad
\partial\overline{S_m}\subset \{|z|=2R_{m-1}\}\cup\{|z|=R_m\}\subset A_{m-1}\cup A_m\ (m\ge2).
\]
Since $g(x)$ avoids every $A_m$, we have
\[
g(x)\cap\partial\overline{S_m}=\varnothing\qquad \mbox{ for all } x\in S^k,\ \mbox{ for all } m\ge1.
\]

\medskip
\noindent\textbf{Step 2 (restrict to each shell as a finite configuration map).}
For each $m$, apply the restriction lemma (Lemma~\ref{lem:restrict}) with
$K=S^k$, $F=g$, and $B=\overline{S_m}$.
Because $S^k$ is connected, we obtain an integer $n_m\ge0$ such that
\[
n_m=\#\bigl(g(x)\cap \overline{S_m}\bigr)\qquad \mbox{ for all } x\in S^k,
\]
and a continuous map
\[
g_m:S^k\to C_{n_m}(S_m),\qquad g_m(x):=g(x)\cap S_m.
\]
Let $Q_m:=g_m(*)=Q\cap S_m$, so that (disjoint union)
\[
Q=\bigsqcup_{m\ge1} Q_m.
\]

\medskip
\noindent\textbf{Step 3 (each shell map is null-homotopic).}
Each $S_m$ is an open disk ($m=1$) or an open annulus ($m\ge2$). By assumption,
each $C_{n_m}(S_m)$ is aspherical. Hence for $k\ge2$,
\[
\pi_k\bigl(C_{n_m}(S_m),Q_m\bigr)=0,
\]
so $g_m:(S^k,*)\to (C_{n_m}(S_m),Q_m)$ admits a based null-homotopy
\[
H_m:(S^k\times[0,1],\{*\}\times[0,1])\to (C_{n_m}(S_m),Q_m)
\]
with $H_m(\cdot,0)=g_m$ and $H_m(\cdot,1)\equiv Q_m$.

\medskip
\noindent\textbf{Step 4 (assemble the shell null-homotopies into a global null-homotopy).}
Let $t_m:=1-2^{-m}$ and $t_0:=0$.
Define $H:S^k\times[0,1]\to C^{lf}_\infty(\C)$ by the standard shell-by-shell concatenation:
for $t\in[t_{m-1},t_m]$, set
\[
H(x,t)\cap S_j=
\begin{cases}
Q_j, & j<m,\\[6pt]
H_m\!\left(x,\dfrac{t-t_{m-1}}{t_m-t_{m-1}}\right), & j=m,\\[10pt]
g_j(x), & j>m,
\end{cases}
\qquad
H(x,t)\cap\Bigl(\C\setminus\bigcup_{j\ge1}S_j\Bigr)=\varnothing.
\]
This is well-defined because the $S_j$ are disjoint.

\smallskip
\emph{(a) $H(x,t)\in C^{lf}_\infty(\C)$ for all $(x,t)$.}
Local finiteness: every compact set $K_0\subset\C$ meets only finitely many shells $S_j$
(since $R_j\to\infty$), and each $H(x,t)\cap S_j$ is finite, hence $H(x,t)\cap K_0$ is finite.

Infiniteness: any infinite locally finite subset of $\C$ is unbounded, hence meets infinitely many shells.
Thus $Q=g(*)$ meets infinitely many shells, i.e.\ $Q_m\neq\varnothing$ for infinitely many $m$.
For each fixed $m$, $n_m=\#(g(x)\cap S_m)$ is constant in $x$, so $n_m=\#Q_m$.
Therefore $n_m>0$ for infinitely many $m$, and for every $(x,t)$ the union over all $j$
contains points in infinitely many shells; hence $H(x,t)$ is infinite.

\smallskip
\emph{(b) $H$ is continuous in the vague topology.}
Fix $\varphi\in C_c(\C)$ and let $C_0:=\supp(\varphi)$.
Since $C_0$ is compact, there exists $M$ such that $C_0\cap S_j=\varnothing$ for all $j>M$.
Hence
\[
\langle H(x,t),\varphi\rangle=\sum_{j=1}^{M}\ \sum_{z\in H(x,t)\cap S_j}\varphi(z),
\]
a finite sum. On each slab $[t_{m-1},t_m]$, only the $m$-th shell term varies with $t$;
all others are either fixed $Q_j$ or fixed $g_j(x)$.
Because each $H_m$ is continuous into the finite configuration space $C_{n_m}(S_m)$,
the functional $A\mapsto \sum_{a\in A}\varphi(a)$ is continuous on $C_{n_m}(S_m)$, and therefore
$(x,t)\mapsto \langle H(x,t),\varphi\rangle$ is continuous on each slab.
At the junction times $t=t_m$, continuity holds since $H_m(\cdot,1)=Q_m$ matches the next slab's value.
Near $t=1$, all shells meeting $C_0$ have already been frozen to $Q_j$, so the expression stabilizes.
Thus $\langle H(\cdot,\cdot),\varphi\rangle$ is continuous; since the vague topology is initial for all
$\varphi\in C_c(\C)$, $H$ is continuous.

\smallskip
We have $H(\cdot,0)=g$, $H(\cdot,1)\equiv Q$, and $H(*,t)=Q$ for all $t$
(based condition because each $H_m$ is based).
Hence $g$ is based null-homotopic at basepoint $Q$, so
\[
[g]=0\in \pi_k\bigl(C^{lf}_\infty(\C),Q\bigr).
\]

\medskip
\noindent\textbf{Step 5 (return to basepoint $\N$).}
Let $F:S^k\times[0,1]\to C^{lf}_\infty(\C)$ be a homotopy from $f$ to $g$, and set $\gamma(t):=F(*,t)$,
a path from $\N$ to $Q$.
By basepoint change (Lemma~\ref{lem:basepoint-change}),
\[
\gamma_\#([f])=[g]=0.
\]
Since $\gamma_\#$ is an isomorphism for $k\ge2$, we conclude $[f]=0$ in
$\pi_k(C^{lf}_\infty(\C), \N)$.
As $f$ was arbitrary, $\pi_k(C^{lf}_\infty(\C), \N)=0$ for all $k\ge2$.
\end{proof}

\section{Asphericity of $Conf^{lf}_{\infty}(\C)$}

\subsection{Parametrized isotopy actions and continuity}

\begin{lemma}[Uniform bound on points in a compact]\label{lem:compact-count}
Let $B_n,B\in C^{lf}_\infty(\C)$ and suppose $B_n\to B$ vaguely.
Then for every compact $K\subset\C$,
\[
\sup_n \#(B_n\cap K)<\infty.
\]
\end{lemma}

\begin{proof}
Choose $\chi\in C_c(\C)$ such that $\chi\ge 1$ on $K$. Then for each $n$,
\[
\#(B_n\cap K)\le \sum_{b\in B_n}\chi(b),
\]
and the right-hand side converges to $\sum_{b\in B}\chi(b)<\infty$ by vague convergence.
Hence it is bounded in $n$.
\end{proof}

Let $K$ be a compact space.
A \emph{$K$--parametrized compactly supported ambient isotopy} of\/ $\C$ means a continuous map
\[
\Phi:K\times [0, 1]\times \C \to \C,\qquad (x,t,z)\mapsto \Phi_{x,t}(z),
\]
such that for each $(x,t)$ the map $\Phi_{x,t}:\C\to\C$ is a homeomorphism, $\Phi_{x,0}=\mathrm{id}$ for all $x$,
and there exists a \emph{single} compact set $K_{\mathrm{iso}}\subset\C$ with
\[
\Phi_{x,t}(z)=z\qquad \text{for all }(x,t)\in K\times[0,1]\text{ and all }z\in \C\setminus K_{\mathrm{iso}}.
\]
Given $B\in C^{lf}_\infty(\C)$ define $\Phi_{x,t}(B):=\{\Phi_{x,t}(b):b\in B\}$.
Given $A=(a_j)_{j\geq 1}\in\Conf$ define $\Phi^\infty_{x,t}(A):=(\Phi_{x,t}(a_1),\Phi_{x,t}(a_2),\dots)$.
Then
\[
P(\Phi_{x,t}^\infty(A))=\Phi_{x,t}(P(A)).
\]

\begin{lemma}[Parametrized isotopy action is continuous on $C^{lf}_\infty(\C)$]\label{lem:isotopy-unlabeled-param}
Let $K$ be compact and let $\Phi$ be a $K$--parametrized compactly supported ambient isotopy as above.
Then the map
\[
C^{lf}_\infty(\C)\times K\times[0,1]\to C^{lf}_\infty(\C),
\qquad (B,x,t)\mapsto \Phi_{x,t}(B),
\]
is continuous in the vague topology.
\end{lemma}

\begin{proof}
Fix $\phi\in C_c(\C)$ and set
\[
\Lambda_\phi(B,x,t):=\sum_{z\in \Phi_{x,t}(B)}\phi(z)=\sum_{b\in B}\phi(\Phi_{x,t}(b)).
\]
It suffices to prove $\Lambda_\phi$ is continuous.

Let $(B_n,x_n,t_n)\to(B,x,t)$ where $B_n\to B$ vaguely and $(x_n,t_n)\to(x,t)$.
Write
\begin{align*}
\Lambda_\phi(B_n,x_n,t_n)-\Lambda_\phi(B,x,t)
&=\Bigl(\sum_{b\in B_n}\phi(\Phi_{x_n,t_n}(b))-\sum_{b\in B_n}\phi(\Phi_{x,t}(b))\Bigr)
+\Bigl(\sum_{b\in B_n}\psi(b)-\sum_{b\in B}\psi(b)\Bigr)\\
&=:(\star)+(\star\star),
\end{align*}
where $\psi:=\phi\circ \Phi_{x,t}\in C_c(\C)$ is fixed. Then $(\star\star)\to 0$ by vague convergence.

For $(\star)$, set the compact $K_0:=(\supp\phi)\cup K_{\mathrm{iso}}$.
If $b\notin K_0$, then $b\notin K_{\mathrm{iso}}$ implies $\Phi_{x_n,t_n}(b)=\Phi_{x,t}(b)=b$,
and also $b\notin\supp\phi$, so both summands vanish. Hence
\[
|(\star)|\le \#(B_n\cap K_0)\cdot \|\phi\circ \Phi_{x_n,t_n}-\phi\circ \Phi_{x,t}\|_{L^\infty(K_0)}.
\]
By Lemma~\ref{lem:compact-count}, $\#(B_n\cap K_0)$ is uniformly bounded in $n$.
Since $\Phi$ is continuous on the compact set $K\times [0, 1]\times K_0$, we have uniform convergence
$\Phi_{x_n,t_n}\to \Phi_{x,t}$ on $K_0$, hence
$\phi\circ \Phi_{x_n,t_n}\to \phi\circ \Phi_{x,t}$ uniformly on $K_0$, and thus $(\star)\to 0$.
Therefore $\Lambda_\phi(B_n,x_n,t_n)\to\Lambda_\phi(B,x,t)$, proving continuity.
\end{proof}

\begin{lemma}[Parametrized isotopy action is continuous on $(\Conf,d_{\sum})$]\label{lem:isotopy-ordered-param}
Let $K$ be compact and let $\Phi$ be a $K$--parametrized compactly supported ambient isotopy.
Then the map
\[
\Conf\times K\times[0,1]\to\Conf,\qquad (A,x,t)\mapsto \Phi_{x,t}^\infty(A),
\]
is continuous for $d_{\sum}$.
\end{lemma}

\begin{proof}
Product part: for each coordinate $j$, the map $(A,x,t)\mapsto \Phi_{x,t}(a_j)$ is continuous
(since $(x,t,z)\mapsto \Phi_{x,t}(z)$ is continuous). Hence $d_{\mathrm{prod}}$-continuity follows.

Vague part: $P(\Phi^\infty_{x,t}(A))=\Phi_{x,t}(P(A))$ and $P:\Conf\to C^{lf}_\infty(\C)$ is continuous. 
Thus the composition
\[
(A,x,t)\mapsto (P(A),x,t)\mapsto \Phi_{x,t}(P(A))
\]
is continuous in the vague topology by Lemma~\ref{lem:isotopy-unlabeled-param}.
Therefore the $d_{\mathcal V}$ term varies continuously, and so does $d_{\sum}$.
\end{proof}

\subsection{Infinite gaps in the ordered space}

\begin{lemma}[Ordered gaps]\label{lem:ordered-gaps}
Let $m\ge2$ and let $f:S^m\to(\Conf,d_{\sum})$ be continuous.
Then there exist radii $0<R_1<R_2<\cdots\to\infty$ and a map $\widehat f:S^m\to\Conf$
homotopic to $f$ in $\Conf$ such that for every $\ell\ge1$ and every $x\in S^m$,
\[
P(\widehat f(x))\cap\{z\in\C:\ R_\ell\le |z|\le 2R_\ell\}=\varnothing.
\]
Moreover, the homotopy $f\simeq \widehat f$ may be chosen as a countable concatenation on intervals
$[t_{\ell-1},t_\ell]$ with $t_\ell\uparrow 1$.
\end{lemma}

\begin{proof}
Apply Lemma~\ref{lem:inf-gap} to $F:=P\circ f:S^m\to C^{lf}_\infty(\C)$.
We obtain radii $R_\ell\to\infty$ and a homotopy $F\simeq G$ which is a countable concatenation of stages.
At stage $\ell$ we have an $S^m$--parametrized compactly supported ambient isotopy $\Phi^{(\ell)}$ supported
in a compact annulus, and these supports are pairwise disjoint.

Lift stage $\ell$ to $\Conf$ by acting coordinatewise:
\[
H^{(\ell)}(x,t):=\bigl(\Phi^{(\ell)}_{x,t}\bigr)^\infty\bigl(A^{(\ell)}(x)\bigr),
\]
where $A^{(\ell)}(x)$ denotes the configuration at the start of stage $\ell$.
By Lemma~\ref{lem:isotopy-ordered-param}, each stage is continuous in $d_{\sum}$.

Reparametrize stage $\ell$ to a subinterval $[t_{\ell-1},t_\ell]$ where $t_\ell=1-2^{-\ell}$,
and concatenate to obtain a homotopy $H:S^m\times[0,1)\to\Conf$.
We now check that $H$ extends continuously to $t=1$, yielding the endpoint $\widehat f$.

\smallskip
\noindent\emph{Product part at $t=1$.}
Fix a coordinate $j$. Since $S^m$ is compact and $x\mapsto f_j(x)$ is continuous, $f_j(S^m)$ is compact,
hence contained in some closed ball $\overline{B}_{M_j}$. Because the annular supports of the stages drift to infinity
(and are pairwise disjoint), only finitely many stage supports intersect $\overline{B}_{M_j}$.
Hence the $j$-th coordinate of $H(x,t)$ changes only during finitely many stages, and thus becomes constant for all sufficiently large $\ell$,
uniformly in $x$. Therefore the $d_{\mathrm{prod}}$ limit as $t\to1^-$ exists and equals the coordinatewise stabilized configuration.

\smallskip
\noindent\emph{Vague part at $t=1$.}
Fix $\phi\in C_c(\C)$. Then $\supp(\phi)$ is compact, hence intersects only finitely many stage supports
(again because the annular supports are pairwise disjoint and escape to infinity).
Thus, for $t$ sufficiently close to $1$, no further stage affects any point contributing to $\sum\phi$,
so $t\mapsto \sum_{z\in P(H(x,t))}\phi(z)$ is eventually constant. This gives continuity at $t=1$ in the vague topology.

Consequently $H$ extends to a continuous homotopy on $[0,1]$ with endpoint $\widehat f$.
Finally, $P(\widehat f)=G$ and $G(x)$ avoids every annulus $\{R_\ell\le |z|\le 2R_\ell\}$, proving the lemma.
\end{proof}

\subsection{Shell decomposition and constant index sets}

Fix $m\ge2$ and a continuous map $f:S^m\to\Conf$.
By Lemma~\ref{lem:ordered-gaps}, after homotopy we may assume there exist radii $0<R_1<R_2<\cdots\to\infty$ such that
\[
P(f(x))\cap\{R_\ell\le |z|\le 2R_\ell\}=\varnothing\qquad(\mbox{ for all } x\in S^m,\ \mbox{ for all } \ell\ge1).
\]
Define shells
\[
S_1:=\{z:|z|<R_1\},\qquad
S_\ell:=\{z:2R_{\ell-1}<|z|<R_\ell\}\quad(\ell\ge2).
\]

\begin{lemma}[Index sets are constant on $S^m$]\label{lem:index-sets}
For each $j\in\N$ there exists a unique $\ell(j)\in\N$ such that $f_j(S^m)\subset S_{\ell(j)}$.
Equivalently,
\[
J_\ell:=\{j\in\N:\ f_j(x)\in S_\ell\}
\]
is independent of $x\in S^m$. Moreover, each $J_\ell$ is finite and infinitely many $J_\ell$ are nonempty.
\end{lemma}

\begin{proof}
Fix $j$. The coordinate map $x\mapsto f_j(x)$ is continuous (product term of $d_{\sum}$).
Because $f_j(x)\in P(f(x))$ and $P(f(x))$ meets no gap annulus, $|f_j(x)|$ never equals any boundary radius $R_\ell$ or $2R_\ell$.
Hence $f_j(x)$ cannot cross from one shell to another, so the shell index is locally constant in $x$.
Since $S^m$ is connected, it is constant; call it $\ell(j)$.

Thus $J_\ell=\{j:\ell(j)=\ell\}$ is independent of $x$.
If $J_\ell$ were infinite then the bounded set $S_\ell$ would contain infinitely many points of $P(f(x))$, contradicting local finiteness.
If only finitely many $J_\ell$ were nonempty, then $P(f(x))$ would be bounded and hence finite by local finiteness,
contradicting that it is infinite. Thus infinitely many $J_\ell$ are nonempty.
\end{proof}

Write $n_\ell:=|J_\ell|$ and list $J_\ell=\{j_{\ell,1}<\cdots<j_{\ell,n_\ell}\}$ (possibly $n_\ell=0$).
Adopt $Conf_0(U)=\{\mathrm{pt}\}$.

\begin{definition}[Shell maps]\label{def:shell-maps}
Define
\[
f^{(\ell)}:S^m\to Conf_{n_\ell}(S_\ell),\qquad
f^{(\ell)}(x):=(f_{j_{\ell,1}}(x),\dots,f_{j_{\ell,n_\ell}}(x)).
\]
\end{definition}

\subsection{Shell null-homotopies relative to the basepoint}

\begin{lemma}[Based shell contraction]\label{lem:shell-based-null}
Fix $m\ge2$ and a basepoint $*\in S^m$. For each $\ell$ there exists a homotopy
\[
H^{(\ell)}:S^m\times[0,1]\to Conf_{n_\ell}(S_\ell)
\]
such that
\[
H^{(\ell)}(\cdot,0)=f^{(\ell)},\qquad
H^{(\ell)}(*,t)=f^{(\ell)}(*)\ \ \mbox{ for all } t,\qquad
H^{(\ell)}(\cdot,1)\equiv f^{(\ell)}(*).
\]
\end{lemma}

\begin{proof}
This is exactly the statement that the based homotopy class of $f^{(\ell)}:(S^m,*)\to(Conf_{n_\ell}(S_\ell),f^{(\ell)}(*))$
is trivial. It follows from Proposition~\ref{prop:finite-aspherical}, since $\pi_m(Conf_{n_\ell}(S_\ell),f^{(\ell)}(*))=0$ for $m\ge2$.
\end{proof}

\subsection{Assembling a global based null-homotopy}

Let $t_\ell:=1-2^{-\ell}$ so $t_\ell\uparrow 1$, and let $s_\ell:[t_{\ell-1},t_\ell]\to[0,1]$ be the affine reparametrization.
Fix a basepoint $*\in S^m$.

For each $\ell$ define the constant shell configuration
\[
c^{(\ell)}:=f^{(\ell)}(*)\in Conf_{n_\ell}(S_\ell),
\]
and choose a based shell homotopy $H^{(\ell)}$ as in Lemma~\ref{lem:shell-based-null}.

Define $\mathcal H:S^m\times[0,1)\to\Conf$ as follows.
For $t\in[t_{\ell-1},t_\ell]$ set, for $1\le i\le n_\ell$,
\[
\mathcal H_{j_{\ell,i}}(x,t):=\bigl(H^{(\ell)}(x,s_\ell(t))\bigr)_i,
\]
and for indices $j=j_{r,i}\in J_r$ with $r\neq \ell$ define
\[
\mathcal H_{j}(x,t):=
\begin{cases}
(c^{(r)})_i & \text{if } r<\ell,\\[2pt]
f_{j}(x) & \text{if } r>\ell.
\end{cases}
\]
Because shells are disjoint and each $H^{(\ell)}$ stays in $Conf_{n_\ell}(S_\ell)$, there are no collisions.
Local finiteness holds since a compact subset of $\C$ meets only finitely many shells and each shell contributes finitely many points.
Thus $\mathcal H(x,t)\in\Conf$ for all $(x,t)$.

Let $c^\infty\in\Conf$ be the configuration defined by
\[
c^\infty_{j_{\ell,i}}:=(c^{(\ell)})_i=f_{j_{\ell,i}}(*)
\qquad(\mbox{ for all } \ell,\ 1\le i\le n_\ell).
\]
Equivalently, $c^\infty=f(*)$.

\begin{lemma}[Continuity and extension]\label{lem:global-cont}
The map $\mathcal H:S^m\times[0,1)\to(\Conf,d_{\sum})$ is continuous and extends continuously to
\[
\mathcal H:S^m\times[0,1]\to(\Conf,d_{\sum})
\]
by setting $\mathcal H(x,1):=c^\infty$.
Moreover $\mathcal H(*,t)=c^\infty$ for all $t\in[0,1]$.
\end{lemma}

\begin{proof}
On each slab $S^m\times[t_{\ell-1},t_\ell]$ only finitely many coordinates vary with $t$ (those in $J_\ell$),
so $\mathcal H$ is continuous in $d_{\mathrm{prod}}$ on $[0,1)$.
For the vague part on $[0,1)$, fix $\phi\in C_c(\C)$ with $\supp(\phi)\subset \overline{B}_R$ and choose $L$
such that $\overline{B}_R$ meets only shells $S_1,\dots,S_L$.
Then $(x,t)\mapsto \sum_{z\in P(\mathcal H(x,t))}\phi(z)$ is a finite sum of continuous contributions from shells $1,\dots,L$,
hence continuous; this implies $d_{\mathcal V}$-continuity, and thus $d_{\sum}$-continuity.

For the extension at $t=1$, the product part follows since each coordinate $j$ becomes constant after time $t_{\ell(j)}$.
For the vague part, the same compact-support argument shows that for $t\ge t_L$ the $\phi$--sum is already constant in $t$,
hence equals its value at $t=1$.

Finally, since each shell homotopy is based, $H^{(\ell)}(*,s)=f^{(\ell)}(*)$, the assembled homotopy satisfies
$\mathcal H(*,t)=f(*)=c^\infty$ for all $t$.
\end{proof}

\subsection{Asphericity}

\begin{theorem}\label{thm:aspherical-conf}
The space $(\Conf,d_{\sum})$ is aspherical.
\end{theorem}

\begin{proof}
Fix $m\ge2$ and a based map $f:(S^m,*)\to(\Conf, \widetilde{\N})$.

Apply Lemma~\ref{lem:ordered-gaps} to obtain a homotopy in $\Conf$ from $f$ to a map $\widehat f$
whose unlabeled configurations have infinitely many gaps.
Let $\gamma(t)$ be the basepoint track of this homotopy, a path from $\widetilde{\N}$ to $\widehat f(*)$.
By Lemma~\ref{lem:basepoint-change}, the induced map
\[
\gamma_{\#}:\pi_m(\Conf, \widetilde{\N})\xrightarrow{\ \cong\ }\pi_m(\Conf,\widehat f(*))
\]
is an isomorphism. Thus it suffices to show $[\widehat f]=0$ in $\pi_m(\Conf,\widehat f(*))$.

Perform the shell decomposition for $\widehat f$ (renaming it $f$ for notational simplicity), and construct the global based homotopy
$\mathcal H$ as in Lemma~\ref{lem:global-cont}.
This homotopy satisfies $\mathcal H(\cdot,0)=\widehat f$, $\mathcal H(*,t)=\widehat f(*)$ for all $t$, and
$\mathcal H(\cdot,1)$ is the constant map at $\widehat f(*)$. Hence $\widehat f$ is based null-homotopic and $[\widehat f]=0$.

Therefore $\gamma_{\#}([f])=[\widehat f]=0$, and since $\gamma_{\#}$ is an isomorphism, $[f]=0$ in $\pi_m(\Conf, \widetilde{\N})$.
\end{proof}

\section{Countably infinite sheeted coverings}
\subsection{Homotopy quotient of ordered configuration space}
Throughout this section, we set $G:=\Aut(\N)$, regarded as a discrete group, and $X=\Conf$ with the sum metric $d_{\sum}$.
Recall that for $\phi\in G$, $\phi$ induces a homeomorphism $\phi_*: X\rightarrow X$ by $\phi_*(x_j)_{j\geq 1}=(x_{\phi(j)})_{j\geq 1}$ and the map $P: X \rightarrow \Clf$ is $G$-invariant.

\begin{definition}[Homotopy quotient (Borel construction)]
Fix a contractible free \emph{right} $G$--CW complex $EG$ such that the quotient map $EG\to BG:=EG/G$ is a principal $G$--bundle. Since $G$ is a discrete group, we may take Milnor's model of $EG$.
Define the \emph{homotopy quotient}
\[
X//G := EG\times_G X := (EG\times X)/G,
\]
where $G$ acts on $EG\times X$ on the right by
\[
(e,x)\cdot \sigma := (e\sigma^{-1},\ \sigma\cdot x).
\]
Denote the quotient map by $q:EG\times X\to EG\times_G X$.
\end{definition}

The diagonal right action on $EG\times X$ is free and properly discontinuous.
In particular, the quotient map $q:EG\times X\to EG\times_G X$ is a principal $G$--bundle.

The following observation is straightforward.
\begin{proposition}\label{prop:barP}
There is a well-defined continuous map
\[
\overline P: X//G \longrightarrow C^{lf}_{\infty}(\C),\qquad
\overline P([e,x]) := P(x),
\]
satisfying $\overline P\circ q = P\circ \mathrm{pr}_X$. In particular, we have a commutative diagram:
$$\xymatrix{ EG\times X  \ar[r]^-{\mathrm{pr}_X} \ar[d]_{q} & X \ar[d]^-P \\
X//G \ar[r]^-{\overline{P}} & \Clf}$$
\end{proposition}

\subsection{Canonical embeddings of countably infinite-sheeted coverings}

\begin{definition}
Let $Y$ be a topological space and let $g:Y\to X//G$ be continuous.
Define the pullback principal $G$--bundle
\[
P_g:=g^*(EG\times X)\longrightarrow Y.
\]
Using the standard action of $G$ on the discrete set $\N$, define the associated bundle
\[
\pi_W:W_g:=P_g\times_G \N \longrightarrow Y.
\]
(Equivalently, $((y, e, x)\cdot\sigma,n)\sim((y, e, x),\sigma\cdot n)$.) We call $\pi_W$ the \emph{countable covering associated to $g$}.
\end{definition}

\begin{remark}
Because $\N$ is discrete, the associated bundle $\pi_W:W_g\to Y$ is a covering map
with fiber $\N$.
\end{remark}

\begin{lemma}\label{lem:eval-well-defined}
There is a well-defined continuous map
\[
\mathrm{ev}_g:W_g\to \C,\qquad \mathrm{ev}_g([(y, e, x), n]) := \bigl(x\bigr)_n,
\]
where $(y, e, x)\in P_g\subset Y\times EG\times X$ and $(x)_n$ denotes the $n$-th coordinate of $x\in X$.
\end{lemma}

\begin{proof}
Well-definedness: in $W_g=P_g\times_G\N$ we have $[(y, e, x), n]=[(y, e, x)\cdot\sigma,\sigma^{-1}\cdot n]$.
Moreover, the right action on $P_g$ is induced from the right action on $EG\times X$, therefore
\[
\mathrm{ev}_g([(y, e, x)\cdot\sigma,\sigma^{-1}\cdot n])
= \bigl((x\cdot\sigma)\bigr)_{\sigma^{-1}(n)}
= \bigl((x_{\sigma(j)})_{j\geq 1}\bigr)_{\sigma^{-1}(n)}
= \bigl(x\bigr)_n,
\]
so $\mathrm{ev}_g$ is well-defined.

Continuity is local on $Y$: over any open set $U\subset Y$ admitting a trivialization $P_g|_U\cong U\times G$,
we have $W_g|_U\cong U\times\N$ and $\mathrm{ev}_g$ becomes $(y, e(y), x(y), n)\mapsto (x(y))_n$, which is continuous since
$\N$ is discrete and each coordinate map on $X$ is continuous.
\end{proof}

\begin{definition}[Embedded total space]
Define
\[
\iota_g:W_g\to Y\times \C,\qquad
\iota_g(w):=\bigl(\pi_W(w),\ \mathrm{ev}_g(w)\bigr).
\]
Let $E_g:=\iota_g(W_g)\subset Y\times\C$.
\end{definition}

\begin{lemma}\label{lem:incidence-closed}
Let $Y$ be Hausdorff and let $F:Y\to C^{lf}_{\infty}(\C)$ be continuous for the vague topology.
Then the incidence set
\[
E_F:=\{(y,z)\in Y\times\C:\ z\in F(y)\}
\]
is closed in $Y\times\C$.
\end{lemma}

\begin{proof}
Let $(y_\alpha,z_\alpha)_{\alpha\geq 1}$ be a sequence in $E_F$ converging to $(y,z)\in Y\times\C$.
Assume for contradiction that $z\notin F(y)$.

Since $F(y)$ is locally finite in $\C$, it is closed and discrete, so there exists $r>0$ such that
$\overline{B(z,r)}\cap F(y)=\varnothing$.
Choose $\varphi\in C_c(\C)$ with $\varphi\equiv 1$ on $\overline{B(z, \frac{1}{2}r)}$ and $\supp\varphi\subset B(z, r)$.
Then $\sum_{a\in F(y)}\varphi(a)=0$.

By vague continuity of $F$,
\[
\sum_{a\in F(y_\alpha)}\varphi(a)\longrightarrow \sum_{a\in F(y)}\varphi(a)=0.
\]
On the other hand, since $z_\alpha\to z$, we have $z_\alpha\in B(z,\frac{1}{2}r)$ for all sufficiently large $\alpha$, hence
$\varphi(z_\alpha)=1$. Because $z_\alpha\in F(y_\alpha)$, we obtain
\[
\sum_{a\in F(y_\alpha)}\varphi(a)\ge 1
\]
eventually, contradicting convergence to $0$. Therefore $z\in F(y)$ and $(y,z)\in E_F$.
Thus $E_F$ is closed.
\end{proof}

\begin{theorem}\label{thm:closed-embedded-cover}
Assume that $Y$ is Hausdorff. Let $g:Y\to X//G$ be continuous and let $\pi_W:W_g\to Y$ be the associated countably infinite-sheeted covering.
Then:
\begin{enumerate}
\item $\iota_g:W_g\to Y\times\C$ is a topological embedding, and $\pi_W$ identifies with
the restriction $p_1|_{E_g}:E_g\to Y$ under the homeomorphism $W_g\cong E_g$.
In particular $p_1|_{E_g}$ is a covering space equivalent to $\pi_W$.
\item Let $F:=\overline P\circ g:Y\to C^{lf}_{\infty}(\C)$. Then
\[
E_g=\{(y,z)\in Y\times\C:\ z\in F(y)\}.
\]
In particular $E_g$ is closed in $Y\times\C$.
\end{enumerate}
\end{theorem}

\begin{proof}
\textbf{(1) Local graph description.}
Fix $y_0\in Y$. Because $P_g\to Y$ is a principal $G$--bundle, there exists an open neighborhood $U\ni y_0$
and a homeomorphism $P_g|_U\cong U\times G$ over $U$.
Consequently,
\[
W_g|_U \cong (U\times G)\times_G \N \cong U\times \N,
\]
where the last identification sends $[(y,\sigma),n]$ to $(y,\sigma(n))$.

Under this identification, $\iota_g|_{W_g|_U}$ becomes
\[
U\times \N \longrightarrow U\times \C,\qquad (y,n)\longmapsto (y,x_n(y)),
\]
where $y\mapsto x(y)\in X$ is the continuous ordered configuration obtained by projecting
a local section of $P_g$ to the $X$--factor. Each slice $U\times\{n\}$ maps homeomorphically onto the graph
\[
\Gamma_n:=\{(y,x_n(y)):\ y\in U\}\subset U\times\C.
\]
Because $x(y)\in X$ has pairwise distinct coordinates, the graphs $\Gamma_n$ are disjoint.
Hence $\iota_g$ restricts to a homeomorphism
\[
W_g|_U \cong \bigsqcup_{n\in\N} (U\times\{n\})
\ \xrightarrow{\ \cong\ }\ \bigsqcup_{n\in\N}\Gamma_n = E_g\cap (U\times\C).
\]
Therefore $\iota_g$ is a local homeomorphism onto its image, hence a topological embedding.
Moreover, $p_1:U\times\C\to U$ restricts to the standard covering map
$\bigsqcup_{n\in\N}\Gamma_n\to U$, which corresponds exactly to $W_g|_U\to U$ under the above
homeomorphism. This proves (1).

\noindent
\textbf{(2) Identification of fibers and closedness.}
Let $F=\overline P\circ g:Y\to C^{lf}_{\infty}(\C)$, continuous by Proposition~\ref{prop:barP}.
We claim that for each $y\in Y$,
\[
E_{g,y}:=\{z\in\C:(y,z)\in E_g\} = F(y).
\]
Choose $p\in P_g$ over $y$. By definition of $\overline P$, the configuration $F(y)=\overline P(g(y))$
equals $P(x(p))$, i.e.\ the underlying set of coordinates of $x(p)\in X$.
On the other hand, $E_{g,y}$ is exactly the set of values $\mathrm{ev}_g([p,n])=(x(p))_n$ as $n$ varies.
Hence $E_{g,y}=F(y)$ and therefore
\[
E_g=\{(y,z)\in Y\times\C:\ z\in F(y)\}.
\]
Since $F$ is continuous and $Y$ is Hausdorff, Lemma~\ref{lem:incidence-closed} implies $E_g$ is closed.
\end{proof}

\begin{remark}
Theorem~\ref{thm:closed-embedded-cover} can be read as: a continuous map
$g:Y\to Conf^{lf}_{\infty}(\C)//\Aut(\N)$ canonically determines a countable covering $W_g\to Y$,
together with a closed embedding of its total space into $Y\times\C$ whose fiber over $y$ is the
underlying locally finite subset of $\C$ obtained by forgetting labels in $g(y)$.
\end{remark}

\section{Realizing Weierstrass entire coverings}

Throughout this section, let $X$ be Hausdorff, paracompact, and locally contractible.
Let $q:E\to X$ be a countably infinite-sheeted covering, and suppose we are given an embedding
$E\hookrightarrow X\times\C$ such that $q=p_1|_E$ (the restriction of the first projection)
and $E$ is \emph{closed} in $X\times\C$.

\subsection{Weierstrass entire families and realizing problem}
\begin{definition}
A map $F:X\times\C\to\C$ is a \emph{Weierstrass entire family} if
\begin{enumerate}
\item $F$ is continuous;
\item for each $x\in X$, $F_x(z):=F(x,z)$ is entire in $z$;
\item the zero locus $Z(F):=\{(x,z)\in X\times \C : F(x,z)=0\}$ projects as a countably infinite-sheeted covering
$p_1:Z(F)\to X$.
\end{enumerate}
We call $p_1: Z(F) \rightarrow X$ the Weierstrass entire covering associated to $F$.
\end{definition}

\begin{definition}[Realization problem]
Given such an embedded covering $q:E\to X$, we ask whether there exists a Weierstrass entire family $F$
such that $Z(F)$ is \emph{equivalent} to $E$ as coverings over $X$.
\end{definition}

\begin{proposition}
Fix $n\geq 2$. Let $p:S^n \rightarrow \RP^n$ be the standard double covering and $q:\R \rightarrow S^1$ be defined by
$$q(\theta)=e^{i\theta}$$
Then the product covering map $p\times q: S^n\times \R \rightarrow \RP^n\times S^1$ is not equivalent to any entire covering.
\end{proposition}

\begin{proof}
Assume that the product covering is equivalent to an entire covering over $\RP^n\times S^1$. Then
there exists an embedding $h=(p\times q, r):S^n\times \R \rightarrow (\RP^n\times S^1)\times \C$. Define
$\phi:S^n \rightarrow \C$ by
$$\phi(x)=r(x, 0)$$
Note that
$$h(x, 0)=([x], 1, \phi(x)) \mbox{ and } h(-x, 0)=([x], 1, \phi(-x))$$
Since $h$ is injective, this implies $\phi(x)\neq \phi(-x)$ for any $x\in S^n$ which contradicts the Borsuk-Ulam theorem.
\end{proof}

\subsection{Local graph trivialization and uniform finiteness}

\begin{lemma}[Local graph decomposition]
\label{lem:local-graphs}
For every $x_0\in X$ there exists an open neighborhood $U\ni x_0$ and continuous maps $s_n:U\to\C$ ($n\in\N$)
such that
\[
E\cap (U\times\C)=\bigsqcup_{n\in\N}\Gamma(s_n),\qquad
\Gamma(s_n):=\{(x,s_n(x)):x\in U\},
\]
and $s_n(x)\neq s_m(x)$ for $n\neq m$ and each fixed $x\in U$.
\end{lemma}

\begin{proof}
Choose an evenly covered neighborhood $U\ni x_0$ for the covering $q$; then
$q^{-1}(U)=\bigsqcup_{n\in\N}E_n$ with each $q|_{E_n}:E_n\to U$ a homeomorphism.
Define $s_n(x):=p_2((q|_{E_n})^{-1}(x))$, where $p_2:X\times\C\to\C$ is the second projection.
Then $E_n=\Gamma(s_n)$ and the union is disjoint.
\end{proof}

\begin{lemma}\label{lem:uniform-finite}
Fix $U$ and $E|_U=\bigsqcup_{n\in\N}\Gamma(s_n)$ as in Lemma~\ref{lem:local-graphs}.
Let $K\subset U$ be compact and $R>0$. Then only finitely many graphs $\Gamma(s_n)$ meet
$K\times \overline{B(0,R)}$. Equivalently, there exists $N=N(K,R)$ such that
for any  $x\in K$ and any $n\ge N$, $|s_n(x)|>R$.
\end{lemma}

\begin{proof}
Assume infinitely many $\Gamma(s_n)$ meet $K\times\overline{B(0,R)}$.
Choose points $(x_n,s_n(x_n))$ with $x_n\in K$, $|s_n(x_n)|\le R$, and with indices $n$ all distinct.
Because $K\times\overline{B(0,R)}$ is compact and $E$ is closed in $X\times\C$, a subsequence converges to
\[
(x_n,s_n(x_n))\to (x_\infty,z_\infty)\in E\cap\bigl(K\times\overline{B(0,R)}\bigr).
\]
Since $q$ is a covering, there exists an open neighborhood $\mathcal{O}\subset E$ of $(x_\infty,z_\infty)$ such that
$q|_{\mathcal{O}}:\mathcal{O}\to V$ is a homeomorphism onto an open set $V\ni x_\infty$.
In particular, $\mathcal{O}$ is contained in a \emph{single} sheet of $q^{-1}(V)$.

As $(x_n,s_n(x_n))\to(x_\infty,z_\infty)$ in $E$, for all sufficiently large $n$ we have
$(x_n,s_n(x_n))\in\mathcal{O}$, hence these points lie in the same sheet over $V$, and therefore in the same graph
$\Gamma(s_{n_0})$ for some fixed index $n_0$. This contradicts the choice of distinct indices. Hence only finitely
many graphs meet $K\times\overline{B(0,R)}$, giving the claimed uniform bound.
\end{proof}

\subsection{Local Weierstrass products}

For $p\in\N\cup\{0\}$, define the Weierstrass elementary factor
\[
\mathcal{E}_p(w):=(1-w)\exp\!\left(w+\frac{w^2}{2}+\cdots+\frac{w^p}{p}\right),
\qquad \mathcal{E}_0(w):=(1-w).
\]

\begin{lemma}\label{lem:logEp}
For $|w|<1$ one has
\[
\log\mathcal{E}_p(w)=-\sum_{k=p+1}^{\infty}\frac{w^k}{k},
\]
hence for $|w|\le \tfrac12$,
\[
\bigl|\log\mathcal{E}_p(w)\bigr|\le 2|w|^{p+1}.
\]
\end{lemma}

\begin{proof}
Use $\log(1-w)=-\sum_{k\ge 1}w^k/k$ for $|w|<1$ and cancel the first $p$ terms against the exponential correction.
For $|w|\le 1/2$ bound the tail geometric series:
$\sum_{k\ge p+1}|w|^k/k\le \sum_{k\ge p+1}|w|^k\le 2|w|^{p+1}$.
\end{proof}

\begin{remark}[Avoiding a common point]
\label{rem:choose-c}
If $E\subset X\times\C$ is closed and $(x_0,c)\notin E$, then there exist neighborhoods $U\ni x_0$ and $B(c,\varepsilon)$
such that $(U\times B(c,\varepsilon))\cap E=\varnothing$. In particular, one can locally choose a constant
$c\in\C$ avoided by all fibers over a neighborhood.
\end{remark}

\begin{proposition}\label{prop:local-realization}
Let $U\subset X$ be open and suppose $E|_U=\bigsqcup_{n\in\N}\Gamma(s_n)$ as in Lemma~\ref{lem:local-graphs}.
Assume there exists a constant $c\in\C$ such that $c\neq s_n(x)$ for all $x\in U$ and all $n$.
Then
\[
F_U(x,z):=\prod_{n=1}^{\infty}\mathcal{E}_n\!\left(\frac{z-c}{s_n(x)-c}\right)
\]
converges locally uniformly on $U\times\C$ to a continuous function, and for each $x\in U$ the function $z\mapsto F_U(x,z)$
is entire with simple zeros exactly at $\{s_n(x)\}_{n\in\N}$. Hence $Z(F_U)|_U=E|_U$.
\end{proposition}

\begin{proof}
Fix a compact $K\subset U$ and a radius $R>0$. By Lemma~\ref{lem:uniform-finite} applied to the compact $K$ and radius
$2R$ in the shifted coordinate $s_n-c$, there exists $N$ such that
\[
\mbox{ for all } x\in K\ \mbox{ and } n\ge N:\quad |s_n(x)-c|>2R.
\]
For $|z-c|\le R$ and $n\ge N$, we have $|w_n(x,z)|\le 1/2$, where $w_n(x,z):=(z-c)/(s_n(x)-c)$.
By Lemma~\ref{lem:logEp},
\[
\bigl|\log\mathcal{E}_n(w_n(x,z))\bigr|\le 2|w_n(x,z)|^{n+1}\le 2\cdot 2^{-(n+1)}.
\]
Thus the logarithm series converges uniformly on $K\times\overline{B(c,R)}$ by the $M$--test, and therefore the product
converges uniformly there. Standard Weierstrass-product arguments yield continuity in $(x,z)$ and holomorphy in $z$.

For fixed $x$, each factor has a simple zero at $z=s_n(x)$ and no other zeros; the locally uniform convergence ensures the
limit function has zeros precisely at those points with the same multiplicities. Hence the zero set is exactly
$\bigsqcup_n\Gamma(s_n)$ over $U$.
\end{proof}

\subsection{Global obstruction and criterion}

\begin{definition}
Let $\mathcal{O}$ be the sheaf on $X$ whose sections over an open $U\subset X$ are continuous functions
$h:U\times\C\to\C$ such that $h(x,\cdot)$ is entire for each $x$. Let $\mathcal{O}^{\times}$ be the subsheaf of nowhere-zero
sections.
\end{definition}

\begin{lemma}
\label{lem:derivative-cont}
If $h:U\times\C\to\C$ is continuous and $h(x,\cdot)$ is holomorphic for each $x\in U$, then $\partial h/\partial z$
exists and is continuous on $U\times\C$.
\end{lemma}

\begin{proof}
Fix $(x_0,z_0)$. Choose $r>0$ so that $\overline{D(z_0,r)}$ is contained in the region under consideration.
By Cauchy's integral formula,
\[
\frac{\partial h}{\partial z}(x,z_0)=\frac{1}{2\pi i}\int_{|\zeta-z_0|=r}\frac{h(x,\zeta)}{(\zeta-z_0)^2}\,d\zeta.
\]
The integrand depends continuously on $(x,\zeta)$ on the compact set $U'\times\{|\zeta-z_0|=r\}$ for $U'$ a small
neighborhood of $x_0$. Hence the integral depends continuously on $x$, proving continuity at $(x_0,z_0)$.
\end{proof}

\begin{lemma}[Exact exponential sequence]
\label{lem:exp-seq}
The sequence of sheaves on $X$
\[
0\longrightarrow 2\pi i\,\Z \longrightarrow \mathcal{O} \xrightarrow{\exp} \mathcal{O}^{\times} \longrightarrow 1
\]
is exact.
\end{lemma}

\begin{proof}
Exactness at $\mathcal{O}$ is pointwise: for each $(x,z)$, $\exp(\ell)=1$ if and only if $\ell\in 2\pi i\,\Z$.
Surjectivity of $\exp$ as a sheaf map is local: given $h\in\mathcal{O}^{\times}(U)$ and $x_0\in U$, choose a contractible
neighborhood $V\subset U$ of $x_0$ (by local contractibility). Fix $z_0\in\C$; since $V$ is contractible, the map
$x\mapsto h(x,z_0)\in\C^*$ admits a continuous logarithm $L_0$ on $V$. By Lemma~\ref{lem:derivative-cont},
$h'/h$ is continuous on $V\times\C$, and we may define
\[
\ell(x,z):=L_0(x)+\int_{z_0}^z \frac{\partial h/\partial \zeta}{h}(x,\zeta)\,d\zeta,
\]
along the straight segment. Then $\ell\in\mathcal{O}(V)$ and $e^{\ell}=h$ on $V\times\C$.
\end{proof}

\begin{lemma}[Fineness]
\label{lem:fine}
$\mathcal{O}$ is a fine (hence acyclic) sheaf. In particular, $H^1(X;\mathcal{O})=0$.
\end{lemma}

\begin{proof}
$\mathcal{O}$ is a sheaf of modules over the sheaf $C^0_X$ of continuous complex-valued functions on $X$
(by multiplication in the $x$-variable). On a paracompact space, partitions of unity exist and yield fineness.
\end{proof}

\begin{definition}\label{def:obstruction}
Choose an open cover $\{U_i\}$ of $X$ by \emph{evenly covered} sets on which one can also choose
$c_i\in\C$ avoided by all fibers (Remark~\ref{rem:choose-c}). On each $U_i$ choose a local Weierstrass realization
$F_i\in \mathcal{O}(U_i)$ with $Z(F_i)=E|_{U_i}$ (Proposition~\ref{prop:local-realization}).

On overlaps $U_{ij}:=U_i\cap U_j$, the ratio $h_{ij}:=F_i/F_j$ extends to an element of $\mathcal{O}^{\times}(U_{ij})$
(because $F_i$ and $F_j$ have the same simple zeros). The family $\{h_{ij}\}$ is a \v{C}ech $1$--cocycle with values in
$\mathcal{O}^{\times}$, hence determines a class $[h]\in H^1(X;\mathcal{O}^{\times})$.
Define
\[
c_1(E):=\delta([h])\in H^2(X;2\pi i\,\Z)\cong H^2(X;\Z),
\]
where $\delta$ is the connecting homomorphism from the long exact sequence of Lemma~\ref{lem:exp-seq}.
\end{definition}

\begin{lemma}[Well-definedness]
$c_1(E)$ is independent of all choices of cover and local realizations $\{F_i\}$.
\end{lemma}

\begin{proof}
Replacing $F_i$ by $u_iF_i$ with $u_i\in\mathcal{O}^{\times}(U_i)$ replaces $h_{ij}$ by
$h'_{ij}=(u_i/u_j)h_{ij}$, i.e.\ changes the cocycle by a coboundary, so $[h]$ is unchanged.
Refining the cover does not change the class in sheaf cohomology. Therefore $c_1(E)=\delta([h])$ is well-defined.
\end{proof}

\begin{theorem}[Global realization criterion]
\label{thm:criterion}
Let $q:E\to X$ be a closed embedded countable covering. Then $E$ is equivalent to a Weierstrass entire covering
if and only if $c_1(E)=0$ in $H^2(X;\Z)$.
\end{theorem}

\begin{proof}
If such a global Weierstrass entire family $F$ exists, take $F_i:=F|_{U_i}$. Then $h_{ij}\equiv 1$, so $[h]=0$ and hence $c_1(E)=\delta([h])=0$.

Conversely, assume $c_1(E)=0$. By exactness of the long exact cohomology sequence associated to
$0\to 2\pi i\Z\to\mathcal{O}\to\mathcal{O}^{\times}\to 1$ and Lemma~\ref{lem:fine}, the map
$H^1(X;\mathcal{O}^{\times})\xrightarrow{\delta} H^2(X;2\pi i\Z)$ is injective (its kernel is the image of
$H^1(X;\mathcal{O})=0$). Thus $c_1(E)=\delta([h])=0$ implies $[h]=0$.

Therefore the cocycle $\{h_{ij}\}$ is a coboundary: there exist $u_i\in \mathcal{O}^{\times}(U_i)$ such that
$h_{ij}=u_i/u_j$ on $U_{ij}$. Set $\widetilde F_i:=F_i/u_i$. Then $\widetilde F_i=\widetilde F_j$ on overlaps, hence
the $\widetilde F_i$ glue to a global $F\in\mathcal{O}(X)$. Since division by nowhere-zero factors does not change zeros,
$Z(F)|_{U_i}=Z(F_i)=E|_{U_i}$, hence $Z(F)=E$ globally.
\end{proof}

\begin{corollary}
If $H^2(X;\Z)=0$, then every closed embedded countable covering
$E\subset X\times\C$ is realizable as $Z(F)$ for some Weierstrass entire family $F$.
\end{corollary}

\begin{remark}
If $F:X\times\C\to\C$ is continuous and each $F_x$ is entire, then for each $k\ge 0$ the coefficient
\[
a_k(x):=\frac{1}{2\pi i}\int_{|\zeta|=r}\frac{F(x,\zeta)}{\zeta^{k+1}}\,d\zeta
\]
is continuous in $x$ (by uniform continuity on compacta and Cauchy's integral formula), so
$F_x(z)=\sum_{k\ge 0}a_k(x)z^k$ with continuous $a_k$ holds automatically.
\end{remark}

\subsection{Countably infinite-sheeted coverings and Weierstrass entire coverings}

In the following, we set $G:=\Aut(\N)$ and consider $X$ to be a path-connected, locally path-connected, and semilocally simply connected
space.

\begin{definition}
Let $Cov_{\infty}(X)$ denote the set of equivalence classes of countably infinite-sheeted
covering maps $p:E\to X$.
\end{definition}

\begin{definition}
Fix a basepoint $x_0\in X$. For a covering $p:E\to X$ with countable fiber, choose a bijection
$\lambda:\N\to p^{-1}(x_0)$. Path lifting gives an action of $\pi_1(X,x_0)$ on $p^{-1}(x_0)$,
hence a homomorphism
\[
\rho_{(E,\lambda)}:\pi_1(X,x_0)\longrightarrow \Aut(\N).
\]
Changing $\lambda$ conjugates $\rho_{(E,\lambda)}$ by an element of $\Aut(\N)$.
\end{definition}

\begin{theorem}\label{thm:cov-mono}
There is a bijection
\[
Cov_{\infty}(X)\ \cong\ \Hom(\pi_1(X,x_0),\Aut(\N))\big/\text{conjugacy}.
\]
\end{theorem}

\begin{proof}
This is the standard classification of coverings by $\pi_1$--actions on the fiber.
Given $p:E\to X$, monodromy yields a conjugacy class of homomorphisms.
Conversely, given a group homomorphism $\rho:\pi_1(X,x_0)\to\Aut(\N)$, form the associated covering
$\widetilde X\times_{\rho}\N\to X$, where $\widetilde X$ is the universal cover.
The two constructions are inverse up to isomorphism.
\end{proof}

\begin{definition}
Define
\[
H^{lf}(\infty):=\pi_1(\Conf, \widetilde{\N}),
\qquad
B^{lf}(\infty):=\pi_1(\Conf//G,[e, \widetilde{\N}])
\] where $e$ is the base point of $EG$.
\end{definition}

\begin{proposition}\label{prop:pi1-exact}
There is a locally finite version of the braid sequence:
\[
1\longrightarrow H^{lf}(\infty)\longrightarrow B^{lf}(\infty)\longrightarrow \Aut(\N)\longrightarrow 1,
\]
where the surjection is induced by $\pi:\Conf//G\to BG$. Furthermore, $\Conf$ is an Eilberg-Mac Lane space of type $K(H^{lf}(\infty), 1)$ and $\Conf//G$ is an Eilenberg-Mac Lane space of type $K(B^{lf}(\infty), 1)$.
\end{proposition}

\begin{proof}
Use the long exact homotopy sequence of the fibration
$$\xymatrix{ \Conf\ar[r] & \Conf//G \ar[d]\\
& BG}$$
Since $BG$ is $K(G,1)$, we have $\pi_k(BG)=0$ for $k\ge2$ and
$\pi_1(BG)=G$. By Theorem \ref{thm:aspherical-conf}, $\Conf$ is an Eilberg-Mac Lane space of type $K(H^{lf}(\infty), 1)$. The long exact homotopy sequence implies that $\Conf//G$ is an Eilenberg-Mac Lane space of type $K(B^{lf}(\infty), 1)$ and also gives the short exact sequence.
\end{proof}

\medskip

\begin{definition}
Let $EmbCov_{\infty}^{lf}(X)$ be the set of isomorphism classes of pairs $(p:E\to X,\iota)$ where
$p$ is a countably infinite-sheeted covering and $\iota:E\hookrightarrow X\times\C$ is an embedding over $X$
(i.e.\ $p=p_1\circ\iota$) such that each fiber $\iota(p^{-1}(x))\subset\C$ is locally finite and $\iota(E)$ is
closed in $X\times \C$.
\end{definition}

\begin{theorem}\label{thm:embcov-borel}
Assume $X$ has the homotopy type of a CW complex.
Then there is a natural bijection
\[
[X,\Conf//G]\ \cong\ EmbCov_{\infty}^{lf}(X),
\]
sending a map $g:X\to\Conf//G$ to the embedded covering $(W_g\to X,\iota_g)$ constructed in
Theorem~\ref{thm:closed-embedded-cover}.
\end{theorem}

\begin{proof}
By Theorem \ref{thm:closed-embedded-cover}, a map $g:X\to \Conf//G$ determines a countably infinite-sheeted covering
$E_g\subset X\times \C$ over $X$ and $E_g$ is closed in $X\times \C$.

Conversely, given an embedded locally finite covering $E\subset X\times\C$ which is closed, let $P\to X$ be the principal
$G$--bundle of \emph{labelings} of fibers (bijections $\N\to p^{-1}(x)$). The embedding into $\C$ turns such
a labeling into an ordered locally finite configuration, i.e.\ gives a $G$--equivariant map $P\to\Conf$.
This yields a classifying map $g:X\to\Conf//G$, and the two constructions are inverse up to isomorphism.
\end{proof}

\begin{definition}[Weierstrass-realizable classes]
Let $\mathscr{E}(X)\subset EmbCov_{\infty}^{lf}(X)$ be the subset of equivalence classes
of pairs $(p:E\to X,\iota)$ such that the embedded covering $p_1: \iota(E) \rightarrow X$ is a  
Weierstrass entire coverings.
\end{definition}

The following is a rephrasing of Theorem~\ref{thm:criterion} and its corollary for the present notation.
\begin{theorem}[Chern obstruction criterion]
\label{thm:weier-criterion}
For $(E\to X,\iota)\in EmbCov_{\infty}^{lf}(X)$, the class lies in $\mathscr{E}(X)$ if and only if
$c_1(E)=0\in H^2(X;\Z)$.
In particular, if $H^2(X;\Z)=0$, then $\mathscr{E}(X)=EmbCov_{\infty}^{lf}(X)$.
\end{theorem}

From Proposition~\ref{prop:pi1-exact}, we have a short exact sequence:
\[
1\longrightarrow H^{lf}(\infty)\xrightarrow{\ \iota\ } B^{lf}(\infty)\xrightarrow{\ \tau\ } G\longrightarrow 1
\]
and from results above, we have a diagram:
\[
\xymatrix{
\Hom(\pi_1(X,x_0),H^{lf}(\infty)) \ar[d]_{\iota_*} &  & \\
\Hom(\pi_1(X,x_0),B^{lf}(\infty))^{conj} \ar[d]_{\tau_*} \ar[r]^-{\Phi}
& EmbCov_{\infty}^{lf}(X) \ar[d]^{\mathrm{For}} &
\mathscr{E}(X) \ar[l] \ar[ld]_{\mathrm{For}|_{\mathscr{E}}}
\\
\Hom(\pi_1(X,x_0),G)^{conj} \ar[r]^-{\Psi}
& Cov_{\infty}(X) &
}
\]
where the maps are defined as follows:
\begin{enumerate}
\item $\iota_*:\Hom(\pi_1(X,x_0),H^{lf}(\infty))\to \Hom(\pi_1(X,x_0),B^{lf}(\infty))^{conj}$
is postcomposition with $\iota$ followed by passing to conjugacy class:
\[
\iota_*(\alpha):=[\,\iota\circ \alpha\,].
\]

\item $\tau_*:\Hom(\pi_1(X,x_0),B^{lf}(\infty))^{conj}\to \Hom(\pi_1(X,x_0),G)^{conj}$
is induced by $\tau$:
\[
\tau_*([\varphi]):=[\,\tau\circ \varphi\,].
\]

\item $\Phi:\Hom(\pi_1(X,x_0),B^{lf}(\infty))^{conj}\to EmbCov_{\infty}^{lf}(X)$
is the map obtained by realizing $[\varphi]$ by a free homotopy class of maps
$g:X\to \Conf//G$ (using that $\Conf//G$ is a $K(B^{lf}(\infty),1)$)
and then sending $g$ to the embedded covering $(W_g\to X,\iota_g)$ of Theorem~\ref{thm:embcov-borel}.
Equivalently: choose $g$ with $g_*=\varphi$ and set $\Phi([\varphi]):=[(W_g,\iota_g)]$.

\item $\mathrm{For}:EmbCov_{\infty}^{lf}(X)\to Cov_{\infty}(X)$ forgets the embedding:
\[
\mathrm{For}([(p:E\to X,\iota)]) := [p:E\to X].
\]

\item $\Psi:\Hom(\pi_1(X,x_0),G)^{conj}\to Cov_{\infty}(X)$ is the standard classification map
(Theorem~\ref{thm:cov-mono}), sending $[\rho]$ to the isomorphism class of the covering
$\widetilde{X}\times_\rho \N\to X$.
\end{enumerate}
The diagonal arrow $\mathscr{E}(X)\to Cov_\infty(X)$ is the restriction of $\mathrm{For}$.

We have the following result. This result is inspired by the diagram in \cite[Page 108]{HansenBraidsCoverings}.
\begin{theorem}\label{thm:diagram-comm}
The above diagram of sets is commutative (both the square and the triangle).
\end{theorem}

\begin{proof}
\textbf{(A) Well-definedness of $\Phi$.}
Because $\Conf//G$ is an Eilenberg--Mac Lane space $K(B^{lf}(\infty),1)$ (Proposition~\ref{prop:pi1-exact}),
the standard $K(\pi,1)$ correspondence gives a bijection
\[
[X,\Conf//G]\ \cong\ \Hom(\pi_1(X,x_0),B^{lf}(\infty))^{conj},
\]
where $[X,\Conf//G]$ denotes free homotopy classes.
Composing the inverse of this bijection with the bijection
$[X,\Conf//G]\cong EmbCov_{\infty}^{lf}(X)$ from Theorem~\ref{thm:embcov-borel}
yields a well-defined map $\Phi$ as stated.

\medskip
\textbf{(B) Commutativity of the square.}
Let $[\varphi]\in \Hom(\pi_1(X,x_0),B^{lf}(\infty))^{conj}$ and choose $g:X\to\Conf//G$
with $g_*=\varphi$. Consider the associated principal $G$--bundle
\[
P_g:=g^*(EG\times \Conf)\longrightarrow X
\]
and the associated covering
\[
W_g:=P_g\times_G \N \longrightarrow X.
\]
Fix $p_0\in (P_g)_{x_0}$ and use it to label the fiber by
\[
\lambda_0:\N \xrightarrow{\ \cong\ } (W_g)_{x_0},\qquad n\longmapsto [p_0,n].
\]
For a loop $\gamma$ at $x_0$, let $\widetilde\gamma$ be the lift of $\gamma$ to $P_g$ with $\widetilde\gamma(0)=p_0$.
Then there exists a unique $\sigma(\gamma)\in G$ with $\widetilde\gamma(1)=p_0\cdot \sigma(\gamma)$, and
the lifted path in $W_g$ starting at $[p_0,n]$ ends at
\[
[p_0\cdot \sigma(\gamma),n]=[p_0,\sigma(\gamma)(n)].
\]
Hence the monodromy of the covering $W_g\to X$ (with respect to $\lambda_0$) is
\[
\rho_{(W_g,\lambda_0)}:\pi_1(X,x_0)\to G,\qquad [\gamma]\mapsto \sigma(\gamma).
\]

On the other hand, the canonical projection $\bar\pi:\Conf//G\to BG$ induces on fundamental groups
the surjection $\bar\pi_*=\tau:B^{lf}(\infty)\to G$ (Proposition~\ref{prop:pi1-exact}),
and the principal bundle $P_g$ is isomorphic to the pullback $(\bar\pi\circ g)^*(EG)\to X$.
Therefore the principal-bundle monodromy equals
\[
(\bar\pi\circ g)_*=\bar\pi_*\circ g_*=\tau\circ \varphi.
\]
By the computation above, this is exactly $\rho_{(W_g,\lambda_0)}$.
Passing to conjugacy classes (to forget the labeling) shows that the isomorphism class of the underlying covering
$\mathrm{For}(\Phi([\varphi]))=[W_g\to X]$ corresponds under Theorem~\ref{thm:cov-mono} to
\[
[\rho_{(W_g,\lambda_0)}]=[\tau\circ \varphi]=\tau_*([\varphi]).
\]
Equivalently,
\[
\mathrm{For}\circ \Phi \;=\; \Psi\circ \tau_*,
\]
which is the commutativity of the square.

\medskip
\textbf{(C) Commutativity of the triangle.}
By definition $\mathscr{E}(X)\subset EmbCov_{\infty}^{lf}(X)$ and the diagonal arrow
$\mathscr{E}(X)\to Cov_\infty(X)$ is $\mathrm{For}$ restricted to $\mathscr{E}(X)$.
Hence the triangle commutes tautologically.
\end{proof}

\begin{corollary}
$EmbCov_{\infty}^{lf}(X)=Cov_{\infty}(X)$ if and only if $\tau_*$ is surjective. This is true in particular if
$\pi_1(X)$ is a free group.
\end{corollary}

\begin{theorem}
If $\pi_1(X)$ is a free group and $H^2(X; \Z)=0$, then $\mathscr{E}(X)=Cov_{\infty}(X)$.
\end{theorem}

\begin{proof}
Let $p:E\to X$ be any countably infinite-sheeted covering. By Theorem~\ref{thm:cov-mono}, it corresponds to a monodromy
homomorphism $\rho:\pi_1(X)\to G$ up to conjugacy. The group $\pi_1(X)$ is free, so we may choose for each
free generator a lift in $B^{lf}(\infty)$ along the surjection
$B^{lf}(\infty)\twoheadrightarrow G$ (Proposition~\ref{prop:pi1-exact}) and thereby obtain a homomorphism
$\widetilde\rho:\pi_1(X)\to B^{lf}(\infty)$ with $\rho$ as its composite to $G$.

Since $\Conf//G$ is a $K(B^{lf}(\infty),1)$, the homomorphism $\widetilde\rho$ is realized by some map
$g:X\to \Conf//G$. The associated covering $W_g\to X$ has monodromy $\rho$, hence is isomorphic to $E\to X$
(as coverings). Moreover, $W_g$ admits the canonical closed locally finite embedding
$W_g\hookrightarrow X\times\C$ (Theorem~\ref{thm:closed-embedded-cover}).
Finally, $H^2(X;\Z)=0$ implies $c_1(W_g)=0$, so $W_g$ is Weierstrass-realizable by
Theorem~\ref{thm:weier-criterion}.
\end{proof}

For polynomial coverings and analogous results, see \cite[Theorem 6.3]{HansenBraidsCoverings}.

\begin{corollary}\label{cor:graph-case}
If $X$ is a connected graph (equivalently a connected $1$--dimensional CW complex), then
\[
\mathscr{E}(X)=Cov_{\infty}(X).
\]
\end{corollary}

\begin{proof}
Since $X$ is a graph, $\pi_1(X)$ is free and $H^2(X;\Z)=0$. The conclusion follows directly from the above result.
\end{proof}

\section{Appendix}
\begin{proposition*}[{\ref{prop:finite-vague-homeo}}]
For any open $U\subset\C$ and $n\ge 1$, there is a natural homeomorphism
\[
\overline F : C_n(U) \xrightarrow{\ \cong\ } M_n(U),
\]
where $C_n(U)= Conf_n(U)/S_n$ has the quotient topology and $M_n(U)$ has the vague topology.
\end{proposition*}

\begin{proof}
Define
\[
F: Conf_n(U)\to M_n(U),\qquad
F(z_1,\dots,z_n):=\sum_{i=1}^n\delta_{z_i}.
\]
This map is $S_n$--invariant, hence descends to a bijection
\[
\overline F:C_n(U)\to M_n(U),
\qquad
\overline F(\{z_1,\dots,z_n\})=\sum_{i=1}^n\delta_{z_i}.
\]

\medskip
\noindent\emph{Step 1: $\overline F$ is continuous.}
Fix $\varphi\in C_c(U)$. The functional
\[
I_\varphi:M_n(U)\to\R,\qquad I_\varphi(\mu)=\int\varphi\,d\mu
\]
is continuous by definition of the vague topology. On $ Conf_n(U)$,
\[
T_\varphi: Conf_n(U)\to\R,\qquad T_\varphi(z_1,\dots,z_n)=\sum_{i=1}^n\varphi(z_i)
\]
is continuous, and $I_\varphi\circ F=T_\varphi$. Hence $F$ is continuous.
Since $F=\overline F\circ q$ with $q: Conf_n(U)\to C_n(U)$ the quotient map, $\overline F$ is continuous.

\medskip
\noindent\emph{Step 2: $\overline F^{-1}$ is continuous.}
Fix
\[
\mu=\sum_{i=1}^n\delta_{z_i}\in M_n(U),
\qquad z_i\neq z_j\ (i\neq j).
\]
Choose pairwise disjoint open discs $D_1,\dots,D_n\subset U$ such that
\[
z_i\in D_i,\qquad \overline{D_i}\subset U,\qquad \overline{D_i}\cap\overline{D_j}=\varnothing\ (i\neq j).
\]
Next choose pairwise disjoint open sets $\widetilde D_i\subset U$ such that
\[
\overline{D_i}\subset \widetilde D_i,\qquad \overline{\widetilde D_i}\subset U,
\qquad \widetilde D_i\cap \widetilde D_j=\varnothing\ (i\neq j).
\]

For each $i$, choose $\phi_i\in C_c(U)$ with
\[
0\le \phi_i\le 1,\qquad \phi_i(z_i)=1,\qquad \supp(\phi_i)\subset D_i.
\]
Then $\int\phi_i\,d\mu=1$. By vague continuity of $\nu\mapsto\int\phi_i\,d\nu$,
there exists a vague neighborhood $\mathcal N(\mu)\subset M_n(U)$ such that for all $\mu'\in\mathcal N(\mu)$,
\[
\int\phi_i\,d\mu'>\frac12\qquad \text{for all }i.
\]
Write $\mu'=\sum_{j=1}^n\delta_{w_j}$ with distinct $w_j$.
Since $\phi_i$ is supported in $D_i$ and takes values in $[0,1]$,
the inequality $\int\phi_i\,d\mu'>\frac12$ implies that $D_i$ contains at least one atom $w_j$.
Because the discs $D_i$ are pairwise disjoint and $\mu'$ has exactly $n$ atoms,
it follows that for every $\mu'\in\mathcal N(\mu)$:
\begin{equation}\label{eq:one-atom-per-disc-correct}
\text{each }D_i\text{ contains exactly one atom, and there are no atoms outside }\bigcup_{i=1}^n D_i.
\end{equation}
Hence for each $i$ there is a uniquely determined point $w_i(\mu')\in D_i$ such that
\[
\mu'(\,D_i\,)=1,\qquad \mu'=\sum_{i=1}^n \delta_{w_i(\mu')}.
\]

Define
\[
H:\mathcal N(\mu)\to  Conf_n(U),\qquad
H(\mu'):=\bigl(w_1(\mu'),\dots,w_n(\mu')\bigr).
\]
To show $H$ is continuous, fix $i$ and choose $h_i\in C_c(U;\C)$ such that
\[
h_i(z)=z\ \text{for all }z\in \overline{D_i},\qquad \supp(h_i)\subset \widetilde D_i.
\]
Then for any $\mu'=\sum_{j=1}^n\delta_{w_j}\in\mathcal N(\mu)$,
\[
\int h_i\,d\mu'=\sum_{j=1}^n h_i(w_j)=h_i\bigl(w_i(\mu')\bigr)=w_i(\mu'),
\]
because exactly one atom lies in $D_i\subset \widetilde D_i$ and all other atoms lie outside $\widetilde D_i$
(the $\widetilde D_i$ are disjoint), while $h_i(z)=z$ on $\overline{D_i}$.
Since $\mu'\mapsto \int h_i\,d\mu'$ is continuous in the vague topology,
each coordinate $\mu'\mapsto w_i(\mu')$ is continuous, hence $H$ is continuous.

Now set
\[
\overline H:=q\circ H:\mathcal N(\mu)\to C_n(U),
\]
where $q: Conf_n(U)\to C_n(U)$ is the quotient map.
Then, by construction,
\[
\overline F\circ \overline H=\mathrm{id}_{\mathcal N(\mu)}.
\]
Thus $\overline H$ is a continuous local inverse to $\overline F$ at $\mu$, so $\overline F^{-1}$ is continuous at $\mu$.
Since $\mu\in M_n(U)$ was arbitrary, $\overline F^{-1}$ is continuous everywhere.

Therefore $\overline F$ is a homeomorphism.
\end{proof}

\end{document}